\documentclass[11pt]{article}

\usepackage[margin=1in]{geometry}
\usepackage{float}
\usepackage{adjustbox}
\usepackage{booktabs}
\usepackage{xcolor}
\usepackage{amsmath}
\usepackage{amssymb}
\usepackage{caption}
\usepackage{amsthm}
\usepackage{subcaption}
\usepackage{mathtools}
\usepackage{comment}
\usepackage{soul}
\usepackage{authblk}
\usepackage{hyperref}

\usepackage[
  backend=biber,
  style=apa,
  natbib=true,
  url=false,
  isbn=false,
  doi=true,
  eprint=false,
  sortcites=true,
]{biblatex}
\usepackage{graphicx}

\theoremstyle{definition}

\newtheorem{definition}{Definition}[section]
\newtheorem{remark}{Remark}[section]
\newtheorem{theorem}{Theorem}[section]
\newtheorem{example}{Example}[section]
\newtheorem{corollary}{Corollary}[section]
\newtheorem{lemma}{Lemma}[section]
\newtheorem{proposition}{Proposition}[section]
\title{A characterization and closure property of Thorin's class}
\author[1,2]{Thomas K. Kloster}
\affil[1]{Department of Economics and Business Economics, Aarhus University}
\affil[2]{Department of Data Science and Analytics, BI Norwegian Business School}

\begin{document}

\maketitle

\begin{abstract}
    Thorin's class coincides with the smallest class of probability laws on $\mathbb{R}$ that contains all gamma laws and their reflections, while being closed under convolution and weak convergence. In particular, all Thorin laws are self-decomposable. This paper proves an alternative characterization of Thorin's class in terms of what we call the hyperbolic Bochner-Bernstein property. This property corresponds to the notion of hyperbolically completely monotone Laplace transforms on the positive half line, but extended to the real-valued setting. Using the new characterization, it is shown that Thorin's class is closed under multiplication with independent random variables admitting hyperbolically completely monotone densities, which partly settles an open question raised by Bondesson and allows for proving that a large class of possibly skewed real-valued distributions is in Thorin's class. As particular examples, we apply the closure property to show that Johnson's S\textsubscript{U} distribution and the Pearson type IV distribution are in Thorin's class for all parameter choices.
\end{abstract}

\section{Introduction}
Thorin's class of probability laws introduced in \textcite{Thorin1978} extends the class of generalized gamma convolutions to the real line. The generalized gamma convolutions were first introduced in \textcite{Thorin1977Pareto,Thorin1977Lognormal} to study infinite divisibility of some classical probability laws and a comprehensive modern treatment of these and Thorin's class can be found in Bondesson's monograph; \textcite{Bondesson92}. Thorin's class can be characterized as the smallest class that contains all gamma laws on $(0,\infty)$ and their reflections, while being closed under convolution and weak convergence; multivariate extensions based on this definition are studied in \textcite{Barndorff-NielsenMaejimaSato2006,Torricelli2026}. In the following, we shall denote Thorin's class by $\mathcal{T}(\mathbb{R})$ and we note that it is closed under weak limits and affine maps, and contains all Gaussian and degenerate laws. Every law in $\mathcal{T}(\mathbb{R})$ is self-decomposable and thus corresponds to the stationary law of a Lévy-driven Ornstein-Uhlenbeck process.

The main result of this paper is a characterization of $\mathcal{T}(\mathbb{R})$ in terms of a property of the characteristic function of the laws in $\mathcal{T}(\mathbb{R})$, which we term the \emph{hyperbolic Bochner-Bernstein property}. This is the real-line analogue of hyperbolically completely monotone Laplace transforms cf. Remark~\ref{rem:HFL_vs_HCM} below, which characterizes the class of generalized gamma convolutions cf. \textcite{Bondesson1990}.
\begin{definition}\label{def:HFL}
    A probability law $\eta$ on $\mathbb{R}$ is said to be hyperbolically Bochner-Bernstein (HBB) with Bochner-Bernstein measure $\mu$ if there is a finite positive measure $\mu$ on $\mathbb{R}\times \mathbb{R}$ such that, for all $u,w>0$, $\iint\lvert e^{\mathrm{i}(u-w)x-uwy}\rvert \mu(dx,dy)<\infty$ and the Fourier transform of $\eta$, $\mathcal{F}_\eta$, satisfies
    \begin{equation}\label{eq:HFL_property}
    \mathcal{F}_\eta (u)\mathcal{F}_\eta (-w)=\int_{\mathbb{R}}\int_\mathbb{R} e^{\mathrm{i}(u-w)x-uw y}\mu (dx,dy).
    \end{equation}
    For a random variable $X$ with an HBB law, we write $X\in\mathcal{H}_{BB}$.
\end{definition}
 The main result of the paper is then the following theorem.
\begin{theorem}\label{thm:Thorin_HFL_equivalent}
    A law belongs to $\mathcal{T}(\mathbb{R})$ if and only if it is hyperbolically Bochner-Bernstein. 
\end{theorem}
The proof of Theorem~\ref{thm:Thorin_HFL_equivalent} relies on similar methods as those in Chapter~5 and 7 of \textcite{Bondesson92}, but require some additional care since we must treat the complex-valued Fourier transform as opposed to the more convenient Laplace transform in the case of positive random variables. The hyperbolic Bochner-Bernstein property can be equivalently stated as existence of a particular coupling $(Z,Y)$, where $Z\overset{d}{=}X$ and $Y$ carries information about the symmetric part of $Z$. In particular, it is shown in Proposition~\ref{prop:coupling_type_G} that $Y$ is a generalized gamma convolution and that the symmetrization of $X$ has the law of the random variable $\sqrt{2Y}N$, where $N$ is standard Gaussian. 

Having established Theorem~\ref{thm:Thorin_HFL_equivalent}, the following corollary follows from Theorem~\ref{thm:HCM_multiplication_closure} below, which partly settles the first open problem posed in Section~9 of \textcite{Bondesson2015}.
\begin{corollary}\label{cor:HCM_closure_Thorin}
    Let $X$ be a random variable with law in $\mathcal{T}(\mathbb{R})$ and let $Y$ be an independent and strictly positive random variable that admits a hyperbolically completely monotone density. Then the laws of $YX$ and $X/Y$ are also in $\mathcal{T}(\mathbb{R})$.
\end{corollary}
The class of probability laws with hyperbolically completely monotone \emph{densities} is a strict subclass of those with hyperbolically completely monotone \emph{Laplace transforms} and the closure property of Corollary~\ref{cor:HCM_closure_Thorin} is well-known in the case of generalized gamma convolutions where the product part of the result remains valid with $Y$ a generalized gamma convolution due to Theorem~1 of \textcite{Bondesson2015}. For symmetric Thorin laws, the product part also holds with multipliers $\sqrt{Y}$ with $Y$ a generalized gamma convolution cf. Theorem~2 of \textcite{Bondesson2015}. Whether a general version holds for all Thorin laws is posed as an open problem in Section~9 of \textcite{Bondesson2015} and Corollary~\ref{cor:HCM_closure_Thorin} then settles this problem for a large class of multipliers. The departure from the symmetric setting requires a different strategy than the Gaussian mixture representations typically used in, e.g., \textcite{BehmeMaejimaMatsuiSakuma2012,Bondesson2015,Torricelli2026}.

As an application of the developed theory, we apply the closure property of Corollary~\ref{cor:HCM_closure_Thorin} to establish that the law of a certain class of exponential functionals of Lévy processes that typically arise in insurance risk modelling is in $\mathcal{T}(\mathbb{R})$. This is the content of Section~\ref{sec:applications} and as particular examples of this result, we show that Johnson's S\textsubscript{U} distribution as introduced in \textcite{Johnson1949} and the Pearson type IV distribution belong to $\mathcal{T}(\mathbb{R})$. Both of these distributions are frequently applied as flexible models for skewed and heavy-tailed real-valued data, but their Thorin class membership has, to the best of our knowledge, never been shown previously. In fact, infinite divisibility or self-decomposability of Johnson's S\textsubscript{U} distribution appears to be a novel claim, whereas self-decomposability of the Pearson type IV distribution can be deduced from known results (see Example~\ref{ex:Pearson}). In Section~\ref{sec:mixture}, it is argued that the Gaussian mixture approach of the symmetric case fails for these distributions, thus necessitating an alternative approach, which the hyperbolic Bochner-Bernstein property provides.

\section{Preliminaries and classes of distributions}
In this section, we briefly recap the various classes of probability laws that will be used throughout and some of their most important properties. Throughout, we fix a probability space $(\Omega,\mathcal{F},\mathbb{P})$, on which all random variables are defined. If a random variable $X$ has law $\eta=\mathbb{P}\circ X^{-1}$, then we write $X\sim \eta$ and if $\mathcal{S}$ denotes a class of probability laws, then we write $X\in\mathcal{S}$ if $\eta\in \mathcal{S}$. For a probability law $\eta$, we denote its Fourier transform by $\mathcal{F}_\eta(u)=\int e^{\mathrm{i}ux}\eta(dx)$ for $u\in\mathbb{R}$ and its Laplace transform by $\mathcal{L}_\eta(u)=\int e^{-ux}\eta (dx)$ for $u\in(0,\infty)$. When $X\sim \eta$, we often use the notation $\varphi_X(u)$ for the characteristic function of $X$, i.e., $\varphi_X(u)=\mathcal{F}_\eta(u)$. If it is obvious from the context that $\varphi_X$ pertains to $X$, we shall sometimes simply use $\varphi (u)$ for $\varphi_X(u)$ to avoid excessive notation.

We shall be concerned with infinitely divisible probability laws, which are all probability laws $\eta$ with Fourier transform of Lévy-Khintchine type \citep[see Chapter~2 of][]{Sato2013}
\begin{equation}\label{eq:levy_khintchine}
\mathcal{F}_\eta(u) = \exp\left( \mathrm{i}u\beta-\frac{\sigma^2u^2}{2}+\int_\mathbb{R}\left( e^{\mathrm{i}ux}-1-\mathrm{i}ux\mathbf{1}_{(-1,1)}(x)\right)\nu(dx) \right),
\end{equation}
where $\nu$ is the Lévy measure and satisfies $\int(1\wedge x^2)\nu(dx)<\infty$. If the Lévy measure is of the form $\nu(dx)=\tfrac{g(x)}{\lvert x\rvert}dx$ where the function $g(x)$ is non-negative for all $x\in\mathbb{R}$, non-decreasing on $(-\infty,0)$, and non-increasing on $(0,\infty)$, then we say that $\eta$ is \emph{self-decomposable}. The self-decomposable distributions are unimodal and correspond to stationary solutions of Lévy-driven Ornstein-Uhlenbeck equations cf. \textcite{JurekVerwaart1983,SatoYamazato1983,Wolfe1982}.

\begin{definition}[Generalized gamma convolution]\label{def:GGC}
    A probability law $\eta$ on $[0,\infty)$ is a generalized gamma convolution (GGC) if its Laplace transform has the form
    \begin{equation}\label{eq:GGC_laplace}
    \mathcal{L}_\eta(u) =  \exp\left(-au -\int_{(0,\infty)} \log\left( 1+\tfrac{u}{y} \right)U(dy) \right), \quad u\geq 0,
    \end{equation}
    where $a\geq 0$ and $U$ is a non-negative measure, termed the Thorin measure, satisfying 
    \[
    \int_{(0,1]}\lvert \log (y)\rvert U(dy)<\infty, \qquad \int_{(1,\infty)}y^{-1}U(dy)<\infty.
    \]
    For a random variable $X$ with a GGC law, we write $X\in \mathcal{T}(\mathbb{R}_+)$.
\end{definition}

\begin{definition}[Hyperbolic complete monotonicity]\label{def:HCM}
    A function $f$ on $(0,\infty)$ is hyperbolically completely monotone (HCM) if for every $u>0$, the function $v\mapsto f(uv)f(u/v)$ is completely monotone as a function of $w=v+v^{-1}\in [2,\infty)$. For a random variable $X$ with a HCM density, we write $X\in \mathcal{H}_{CM}$.     
\end{definition}
The GGC class was introduced in \textcite{Thorin1977Pareto,Thorin1977Lognormal} and can equivalently be characterized as the closure under weak convergence of the set of all finite convolutions of gamma laws. It is closed under convolution and scaling -- hence affine maps with positive coefficients -- and taking $U\equiv 0$ in \eqref{eq:GGC_laplace} shows that every Dirac distribution $\delta_a$ with $a\geq 0$ is a GGC. Every GGC is infinitely divisible and self-decomposable. The term ``hyperbolic complete monotonicity'' was first labelled in \textcite{Bondesson92}, although the property has been implicitly used in previous works on GGC's, including those of Thorin. Every HCM density corresponds to a GGC law, but $\mathcal{H}_{CM}$ is a strict subset of the non-degenerate laws in $\mathcal{T}(\mathbb{R}_+)$. The GGC laws can, however, be characterized as the laws that admit HCM Laplace transforms; this was in fact the route taken in \textcite{Thorin1977Lognormal} to show that the lognormal distribution is infinitely divisible. Furthermore, it is shown in \textcite{Bondesson92,SteuteHarn2004} that if $X\in \mathcal{T}(\mathbb{R}_+)$ and $Y\in \mathcal{H}_{CM}$ are independent, then $XY\in \mathcal{T}(\mathbb{R}_+)$, which is upgraded in Theorem~1 of \textcite{Bondesson2015} to the fact that for independent $X,Y\in \mathcal{T}(\mathbb{R}_+)$ then $XY\in \mathcal{T}(\mathbb{R}_+)$. 

\begin{definition}[Thorin's class $\mathcal{T}(\mathbb{R})$]\label{def:thorin_class}
    A probability law $\eta$ on $\mathbb{R}$ is in Thorin's class, $\mathcal{T}(\mathbb{R})$, if its Fourier transform has the form
    \begin{equation}\label{eq:thorin_class_def}
    \mathcal{F}_\eta(u) = \exp\left( \mathrm{i}u\beta -\frac{\sigma^2u^2}{2} + \int_{\mathbb{R}\setminus \lbrace 0 \rbrace}\left[\log\left( \frac{x}{x-\mathrm{i}u}\right)-\frac{\mathrm{i}ux}{1+x^2}\right] U(dx)\right),
    \end{equation}
    where $\beta\in\mathbb{R}$, $\sigma\geq0$, and $U$ is a non-negative measure on $\mathbb{R}\setminus \lbrace 0\rbrace$, termed the Thorin measure, such that
    \begin{equation}\label{eq:Thorin_measure_integrability}
    \int_\mathbb{R}\frac{1}{1+x^2}U(dx)<\infty, \quad \text{and} \quad \int_{\lbrace\lvert x\rvert \leq 1 \rbrace}\lvert \log(x^2)\rvert U(dx)<\infty.
    \end{equation}
\end{definition}
The equivalence of Definition~\ref{def:thorin_class} with the weak closure of finite gamma convolutions and their reflections is proven in \textcite{Thorin1978,Bondesson92}. Comparing the representation \eqref{eq:thorin_class_def} to the classical Lévy-Khintchine representation \eqref{eq:levy_khintchine}, we see that $\mathrm{i}u\beta$ represents the drift plus truncation, $\tfrac{\sigma^2 u^2}{2}$ the Gaussian part, and the integral against $U$ the jump structure. It is shown in \textcite{Bondesson92} that any $X\in\mathcal{T}(\mathbb{R})$ is self-decomposable. The Thorin measure $U$ can be decomposed as $\int_\mathbb{R}f(x)U(dx)=\int_0^\infty f(x)U_+(dx)+\int_0^\infty f(-x)U_-(dx)$, where $f$ is any positive measurable function and
\[
    U_+(B) = U(B),\quad U_-(B)=U(-B),\quad B\subseteq (0,\infty).
\]
This decomposition is convenient to describe the positive and negative jumps of processes with laws in $\mathcal{T}(\mathbb{R})$. We note that the class $\mathcal{T}(\mathbb{R}_+)$ in Definition~\ref{def:GGC} is recovered from $\mathcal{T}(\mathbb{R})$ in Definition~\ref{def:thorin_class} by setting $\sigma=0$ and $a=\beta -\int\tfrac{x}{1+ x^2}U(dx)\geq 0$ with $U$ supported on $(0,\infty)$ and under the additional finite variation criterion $\int_{(1,\infty)}y^{-1}U(dy)<\infty$. Conversely, any random variable in $\mathcal{T}(\mathbb{R})$ supported on $[0,\infty)$ is necessarily in $\mathcal{T}(\mathbb{R}_+)$, since an infinitely divisible law on $[0,\infty)$ has $\sigma=0$, non-negative drift, and no negative jumps with $\int_{(0,1]}x\nu (dx)<\infty$ \citep[see][Theorem~24.11]{Sato2013}. Furthermore, as shown in Chapter~7 of \textcite{Bondesson92}, the Lévy measure of \eqref{eq:thorin_class_def} has density $\lvert x\rvert^{-1}\int_0^\infty e^{-\lvert x\rvert y}U_\pm(dy)$, which forces $U_-=0$ and $\int (1-e^{-y})y^{-1}U(dy)<\infty$, which is equivalent to $\int_{(1,\infty)}y^{-1}U(dy)<\infty$ due to \eqref{eq:Thorin_measure_integrability} and therefore the representation \eqref{eq:thorin_class_def} reduces to \eqref{eq:GGC_laplace}.

\section{The hyperbolic Bochner-Bernstein property}\label{sec:HBB}
The primary goal of this section is to establish useful properties of the hyperbolic Bochner-Bernstein property, which carry over to $\mathcal{T}(\mathbb{R})$ due to Theorem~\ref{thm:Thorin_HFL_equivalent}. Corollary~\ref{cor:HCM_closure_Thorin} follows immediately from Theorem~\ref{thm:Thorin_HFL_equivalent} combined with Theorem~\ref{thm:HCM_multiplication_closure} and the proof of Theorem~\ref{thm:Thorin_HFL_equivalent} and Theorem~\ref{thm:HCM_multiplication_closure} are given in Section~\ref{sec:proofs}.

The hyperbolic Bochner-Bernstein property of Definition~\ref{def:HFL} has, to the best of our knowledge, not been explicitly characterized before. We give an intuitive characterization in terms of a bivariate coupling in Proposition~\ref{prop:Thorin_coupling} and some properties of this coupling are illustrated in Examples~\ref{ex:gamma_HFL}, \ref{ex:laplace_coupling}, \ref{ex:normal_coupling} and Proposition~\ref{prop:coupling_type_G}. Remark~\ref{rem:HFL_vs_HCM} shows that we recover hyperbolically completely monotone Laplace transforms in the case of positive random variables. 

\begin{example}\label{ex:gamma_HFL}
    Let $G\sim \operatorname{Gamma}(\alpha,1)$ be a gamma distributed random variable and set $X=\varepsilon\frac{G}{\lambda}$, where $\lambda>0$ and $\varepsilon\in \lbrace -1,1\rbrace$. Then the characteristic function of $X$ is $\varphi(u)=(1-\tfrac{\mathrm{i}\varepsilon u}{\lambda})^{-\alpha}$ and
    \[
    \varphi(u)\varphi(-w) = \left( 1-\frac{\mathrm{i}\varepsilon(u-w)}{\lambda}+\frac{uw}{\lambda^2}\right)^{-\alpha}=\mathbb{E}\left[ e^{\mathrm{i}(u-w)\tfrac{\varepsilon G}{\lambda}-uw\tfrac{G}{\lambda^2}} \right], \quad u,w>0.
    \]
    This corresponds to the joint Fourier-Laplace transform of the bivariate random variable $(\tfrac{\varepsilon G}{\lambda},\tfrac{G}{\lambda^2})$ and thus $X\in \mathcal{H}_{BB}$. Note that the first coordinate has the same distribution as $X$, which is a general fact cf. Proposition~\ref{prop:Thorin_coupling} below. It is then easy to see that if the HBB property preserves convolutions and weak limits, then any $X\in\mathcal{T}(\mathbb{R})$ must have an HBB law by definition. 
\end{example}

\begin{lemma}\label{lem:Fourier_Laplace_unique}
    Let $\eta$ denote a probability law on $\mathbb{R}\times [0,\infty)$ and define the joint Fourier-Laplace transform as 
    \[
    \mathcal{FL}_\eta(p,q) = \int_\mathbb{R}\int_{[0,\infty)}e^{\mathrm{i}px-qy}\eta (dx,dy).
    \]
    Then $\mathcal{FL}_\eta$ uniquely determines $\eta$, i.e., for two laws $\eta_1,\eta_2$ with $\mathcal{FL}_{\eta_1}(p,q)=\mathcal{FL}_{\eta_2}(p,q)$ for all $(p,q)\in\mathbb{R}\times (0,\infty)$, then $\eta_1=\eta_2$.
\end{lemma}
\begin{proof}
    Let $\eta_1,\eta_2$ be two probability measures with transform $\mathcal{FL}_{\eta_1}(p,q)=\mathcal{FL}_{\eta_2}(p,q)$ for all $p\in\mathbb{R}$ and $q>0$. Define $\nu_j^q(dx)=\int e^{-qy}\eta_j(dx,dy)$. Then $\nu_1^q$ and $\nu_2^q$ have the same Fourier transform $\mathcal{FL}_{\eta_1}(\cdot, q)$ and thus the two measures coincide such that, for every bounded and positive function $g$ on $\mathbb{R}$, it holds that
    \[
        \int_\mathbb{R}\int_{[0,\infty)}g(x)e^{-qy}\eta_1(dx,dy)=\int_\mathbb{R}\int_{[0,\infty)}g(x)e^{-qy}\eta_2(dx,dy),
    \]
    for all $q>0$. For such $g$, the two sides are Laplace transforms of the finite measures $\kappa_j^g(B)=\iint g(x)\mathbf{1}_B(y)\eta_j(dx,dy)$ on $[0,\infty)$. Finite measures on $[0,\infty)$ are determined by their Laplace transform and therefore $\kappa_1^g=\kappa_2^g$. Taking $g(x)=\mathbf{1}_A(x)$ for some Borel set $A$ in $\mathbb{R}$, we find that $\eta_1(A\times B)=\eta_2(A\times B)$ and since rectangles generate the product $\sigma$-algebra, we find that $\eta_1$ and $\eta_2$ coincide.  
\end{proof}
\begin{proposition}\label{prop:Thorin_coupling}
    Let $X$ be a random variable with characteristic function $\varphi$. Then the following statements are equivalent:
    \begin{enumerate}
        \item $X\in\mathcal{H}_{BB}$ with Bochner-Bernstein measure $\mu$.
        \item There exists a random vector $(Z,Y)$ with $Y$ almost surely non-negative, such that, for all $u,w>0$,
        \[
        \mathbb{E}\left[ e^{\mathrm{i}(u-w)Z-uwY} \right] = \varphi(u)\overline{\varphi(w)}.
        \]
        \item There exists a random vector $(Z,Y)$ with $Y$ almost surely non-negative such that, for all $u,w>0$, $uX-wX'\overset{d}{=}(u-w)Z+\sqrt{2uwY}N$, with $N$ standard Gaussian and independent of $(Z,Y)$ and $X'$ an independent copy of $X$.  
    \end{enumerate}
    Moreover, the Bochner-Bernstein measure $\mu$ is unique and a probability measure supported on $\mathbb{R}\times [0,\infty)$ whose first marginal is the law of $X$, i.e. $Z\overset{d}{=}X$.  
\end{proposition}
\begin{proof}
    We first show that the first two points are equivalent. To this end, suppose that $X\in\mathcal{H}_{BB}$ with Bochner-Bernstein measure $\mu$ and set $u=w=\sqrt{q}$ for some $q>0$. It follows from \eqref{eq:HFL_property} that 
    \begin{equation}\label{eq:mu_integral_bound}
    \varphi(\sqrt{q})\overline{\varphi(\sqrt{q})} = \int_\mathbb{R}\int_\mathbb{R}e^{-qy}\mu (dx,dy) = \lvert \varphi (\sqrt{q})\rvert^2 \leq 1.
    \end{equation}
    Note, however, that for $\varepsilon>0$ the integrand is non-negative everywhere and bounded below by $e^{q\varepsilon}$ on $(-\infty,-\varepsilon]$ such that \eqref{eq:mu_integral_bound} yields
    \[
    e^{q\varepsilon}\mu(\mathbb{R}\times (-\infty,-\varepsilon] ) \leq \int_\mathbb{R}\int_\mathbb{R}e^{-qy}\mu (dx,dy) \leq 1.
    \]
    In other words, $\mu(\mathbb{R}\times (-\infty,-\varepsilon])\leq e^{-q\varepsilon}$ and sending $q\to\infty$, we conclude that $\mu(\mathbb{R}\times (-\infty,0))=0$ and hence that $Y$ is almost surely non-negative. It also follows that $\mu$ is necessarily a probability measure since the characteristic function $\varphi$ is continuous in 0, so
    \[
    \mu(\mathbb{R}\times [0,\infty)) = \lim_{q\searrow 0}\lvert \varphi (\sqrt{q})\rvert^2 = \lvert \varphi (0)\rvert^2 = 1.
    \]
    We now show that the first marginal of $\mu$ is the law of $X$. To this end, fix $u>0$ and let $w\to 0$ in \eqref{eq:HFL_property}. By dominated convergence this yields
    \[
    \int_{\mathbb{R}}\int_{\mathbb{R}}e^{\mathrm{i}ux}\mu (dx,dy ) = \varphi (u) \lim_{w\to 0}\varphi(-w) = \varphi(u).
    \]
    Hence the first marginal of $\mu$, denoted $\mu^{(1)}$, agrees with $\varphi $ on $(0,\infty)$ and since both sides are characteristic functions (since $\mu$ is a probability measure), conjugate symmetry yields that $\mu^{(1)}$ is the law of $X$. Now let $(Z,Y)$ be a random vector with law $\mu$. Then $Y\geq 0$ almost surely since the second coordinate of $\mu$ is positive and $X\overset{d}{=}Z$ since $\mu^{(1)}$ is the law of $X$. From \eqref{eq:HFL_property}, it holds by assumption that
    \[
    \mathbb{E}\left[ e^{\mathrm{i}(u-w)Z-uwY} \right] = \varphi_Z(u)\varphi_Z(-w) = \varphi(u)\overline{\varphi(w)},
    \]
    such that $(Z,Y)$ satisfies the claimed properties. Conversely, suppose that $(Z,Y)$ is a random vector satisfying point 2 above. Then we can let $\mu$ be the law of $(Z,Y)$ and this immediately yields the HBB property for $X$. It remains to argue that $\mu$ is unique. This follows since $(u,w)\mapsto (u-w,uw)=:(p,q)$ is a bijection of $(0,\infty)\times (0,\infty)$ into $\mathbb{R}\times (0,\infty)$ and therefore \eqref{eq:HFL_property} determines the joint Fourier-Laplace transform $\mathcal{FL}_\mu(p,q)=\iint e^{\mathrm{i}px-qy}\mu(dx,dy)$ for all $p\in\mathbb{R}$ and $q>0$ cf. Lemma~\ref{lem:Fourier_Laplace_unique}.

    We now show that the third claim is equivalent to the second property, which only requires the identity
    \begin{equation}\label{eq:condition_3}
    \mathbb{E}\left[ e^{\mathrm{i}h(u-w)Z+\mathrm{i}h\sqrt{2uwY}N} \right] = \mathbb{E}\left[ e^{\mathrm{i}h(u-w)Z}\mathbb{E}\left[ e^{\mathrm{i}h\sqrt{2uwY}N}\mid(Z,Y)\right] \right] = \mathbb{E}\left[ e^{\mathrm{i}h(u-w)Z-h^2uwY} \right].
    \end{equation}
    Assuming that the third property holds, simply set $h=1$ in \eqref{eq:condition_3} to obtain the second property. Assuming that the second property holds, the right hand side of \eqref{eq:condition_3} with $h>0$ is $\varphi(hu)\overline{\varphi(hw)}$, which is the characteristic function of $uX-wX'$. Conjugate symmetry and the fact that both sides are $1$ in $h=0$ then yields that $(u-w)Z+\sqrt{2uwY}N$ and $uX-wX'$ have the same characteristic function and hence are equal in law.
\end{proof}

\begin{remark}\label{rem:HFL_vs_HCM}
    Suppose that $X\geq 0$ almost surely with Laplace transform $\mathcal{L}(u)=\varphi(\mathrm{i}u)$ for $u>0$. Then one may show that $\mathcal{L}$ is hyperbolically completely monotone when $X\in \mathcal{H}_{BB}$. From Theorem~6.1.1 in \textcite{Bondesson92}, this means $X\in\mathcal{T}(\mathbb{R}_+)$. This can be proved by analytic continuation of \eqref{eq:HFL_property} one variable at a time, from real arguments to purely imaginary arguments. Applying the Schwarz reflection principle in each variable shows that $\varphi_X(u)\varphi_X(-w)-\iint e^{\mathrm{i}(u-w)x-uwy}\mu(dx,dy)$ is identically zero such that $\mathcal{L}(u)\mathcal{L}(w)=\iint e^{-(u+w)x-uwy}\mu(dx,dy)$. Substituting $u=cv$ and $w=\tfrac{c}{v}$ for some $c>0$ the integral is then $\iint e^{-cx(v+v^{-1})}e^{-c^2 y}\mu(dx,dy)$ which is a positive mixture of $e^{-cx(v+v^{-1})}$ ($x$ is positive since $X\geq 0$ almost surely) and therefore completely monotone in $v+v^{-1}$ for $v>0$. By Definition~\ref{def:HCM}, $\mathcal{L}$ is then hyperbolically completely monotone.
\end{remark}

\begin{example}\label{ex:laplace_coupling}
    Consider the case where $X$ has a Laplace distribution, e.g., $X=E_+-E_-$, where $E_{\pm}\sim \operatorname{Gamma}(1,1)$ are independent. Then $\varphi_X(u)=(1+u^2)^{-1}$ and it is easy to see that
    \[
    \varphi_X(u)\overline{\varphi_X(w)} = \frac{1}{(1+u^2)(1+w^2)} = \mathbb{E}\left[ e^{\mathrm{i}(u-w)Z-uwY} \right], \quad u,w>0,
    \]
    where $Z\overset{d}{=}X$ and $Y\overset{d}{=}E_++E_-\sim \operatorname{Gamma}(2,1)$. The coupling $(Z,Y)$ becomes explicitly given as 
    \[
    Y\sim \operatorname{Gamma}(2,1),\quad Z\mid Y=y \sim \operatorname{Uniform}(-y,y).
    \]
    Note that, unlike the positive case in Example~\ref{ex:gamma_HFL}, $Y$ is not a function of $Z$.
\end{example}
\begin{example}\label{ex:normal_coupling}
    Consider the case where $X\sim N(m,\sigma^2)$ is Gaussian. Then the coupling of Proposition~\ref{prop:Thorin_coupling} becomes $Z\overset{d}{=}X$ and $Y\equiv \sigma^2>0$ is degenerate, since
    \[
    \varphi_X(u)\overline{\varphi_X(w)} = e^{\mathrm{i}m(u-w)-\frac{\sigma^2(u-w)^2}{2}}e^{-uw\sigma^2} = \mathbb{E}\left[ e^{\mathrm{i}(u-w)Z-uw\sigma^2} \right].
    \]
\end{example}

\begin{remark}\label{rem:HFL_property}
    The HBB property of Definition~\ref{def:HFL} says that the function constructed from $\mathcal{F}_\eta$ is the joint Fourier-Laplace transform of a measure $\mu$ and Proposition~\ref{prop:Thorin_coupling} says that $\mu$ is a probability measure supported on $\mathbb{R}\times [0,\infty)$ and with the first marginal of $\mu$ corresponding to $\mathcal{F}_\eta$. Define $\mathcal{FL}_\mu(z,q)=\int_{\mathbb{R}\times [0,\infty)}e^{\mathrm{i}zx-qy}\mu (dx,dy)$ for $z\in\mathbb{R}$ and $q\geq 0$ and note that $\mathcal{FL}_\mu$ completely determines $\mu$, cf. Lemma~\ref{lem:Fourier_Laplace_unique}. Note also that $z=u-w$ and $q=uw$ defines a bijection between
    \[
    (u,w)\in (0,\infty)\times (0,\infty) \quad \text{and} \quad (z,q)\in \mathbb{R}\times (0,\infty),
    \]
    via the reparametrization
    $
    u=\frac{z+\sqrt{z^2+4q}}{2}$ and $w = \frac{-z+\sqrt{z^2+4q}}{2}.
    $
    So the HBB property \eqref{eq:HFL_property} is equivalent to 
    \begin{equation}\label{eq:HBB_joint_FL_characterization}
    \mathcal{FL}_\mu(z,q)=\mathcal{F}_\eta\left( \frac{z+\sqrt{z^2+4q}}{2}\right)\overline{\mathcal{F}_\eta\left( \frac{-z+\sqrt{z^2+4q}}{2}\right)}.
    \end{equation}
    The terminology ``Bochner-Bernstein'' arises from the joint Fourier-Laplace transform $\mathcal{FL}_\mu$, which is an exponential mixture of Bochner type in the first coordinate and Bernstein type in the second. 
\end{remark}

\begin{lemma}\label{lem:HFL_closure}
    Let $X,Y\in \mathcal{H}_{BB}$ be independent with respective Bochner-Bernstein measures $\mu_X,\mu_Y$. Then $X+Y\in \mathcal{H}_{BB}$ with Bochner-Bernstein measure $\mu_X\ast\mu_Y$. If $(X_n)_{n\in \mathbb{N}}$ with $X_n\xrightarrow[n\to\infty]{d}X$ with each $X_n\in \mathcal{H}_{BB}$ with Bochner-Bernstein measure $\mu_n$, then $X\in \mathcal{H}_{BB}$ and has Bochner-Bernstein measure $\mu=\lim_{n\to\infty}\mu_n$.
\end{lemma}
\begin{proof}
    Let $\varphi_X,\varphi_Y$ be the characteristic functions of $X,Y\in \mathcal{H}_{BB}$ with associated probability measure $\mu_X,\mu_Y$ in the representation \eqref{eq:HFL_property}. Defining $f(u)=\varphi_X(u)\varphi_Y(u)$, the stated property follows exactly from independence of $X$ and Y and \eqref{eq:HFL_property} 
    \[
    \begin{aligned}
    f(u)f(-w) &= \int_{\mathbb{R}^2}\int_{\mathbb{R}^2}e^{\mathrm{i}(u-w)(x_1+x_2)-uw(y_1+y_2)}\mu_X(dx_1,dy_1)\mu_Y(dx_2,dy_2) \\
    &= \int_{\mathbb{R}^2}e^{\mathrm{i}(u-w)x-uwy}(\mu_X\ast \mu_Y)(dx,dy).
    \end{aligned}
    \]
    For the second property, let $\varphi_n$ denote the characteristic function of $X_n$ and let $\varphi_n\to\varphi$ pointwise, where $\varphi$ is a characteristic function and each $\varphi_n$ satisfies \eqref{eq:HFL_property} with measures $\mu_n$. Define $\mu_n^{(1)}(dx)=\mu_n(dx,\mathbb{R})$ and $\mu_n^{(2)}(dy)=\mu_n(\mathbb{R}\times dy)$ as the marginal distributions of $\mu_n$. It follows from Proposition~\ref{prop:Thorin_coupling} that $\mu_{n}(\mathbb{R}\times [0,\infty))=1$, i.e., that the second coordinate of $\mu_n$ is non-negative. Note then that
    \begin{equation}\label{eq:mu_1_marginal}
    \begin{aligned}
    \int_{\mathbb{R}}e^{\mathrm{i}ux}\mu_n^{(1)}(dx)&=\int_{\mathbb{R}^2}e^{\mathrm{i}ux}\mu_n(dx,dy)\\
    &= \lim_{w\to 0}\int_{\mathbb{R}^2}e^{\mathrm{i}(u-w)x-uwy}\mu_n(dx,dy) \\
    &= \lim_{w\to 0}\varphi_n(u)\varphi_n(-w)\\
    &= \varphi_n(u).
    \end{aligned}
    \end{equation}
    Thus $\mu_n^{(1)}$ is the law of $X_n$ and since $X_n\to X$ in law, it holds that the sequence $(\mu_n^{(1)})_{n\in\mathbb{N}}$ is tight. For the second marginal, let $(A_n,B_n)$ be a bivariate random variable with law $\mu_n$. Then $\mathbb{E}[e^{-qB_n}]=\lvert \varphi_n(\sqrt{q})\rvert^2$ and it holds that $1-\mathbb{E}[e^{-qB_n}]\geq (1-e^{-qM})\mathbb{P}(B_n>M)$ and therefore that
    \[
    \mathbb{P}\left( B_n>M \right)\leq \frac{1-\lvert\varphi_n(\sqrt{q})\rvert^2}{1-e^{-qM}}.
    \]
    Now, because $\varphi_n(\sqrt{q})\to\varphi(\sqrt{q})$ and $\lvert \varphi(\sqrt{q})\rvert^2\to 1$ as $q\to 0$, it follows that $(\mu_n^{(2)})_{n\in\mathbb{N}}$ is uniformly tight. Consequently, $(\mu_n)_{n\in\mathbb{N}}$ is tight on $\mathbb{R}\times [0,\infty)$ and Prokhorov's theorem yields that there is a subsequence such that $\mu_{n_k}\to \mu$ as $k\to\infty$ for some probability measure $\mu$ on $\mathbb{R}\times [0,\infty)$. The function $(x,y)\mapsto e^{\mathrm{i}(u-w)x-uwy}$ is bounded and continuous on $\mathbb{R}\times [0,\infty)$, so
    \[
    \begin{aligned}
    \int_{\mathbb{R}\times [0,\infty)}e^{\mathrm{i}(u-w)x-uwy}\mu(dx,dy) &= \lim_{k\to\infty}\int_{\mathbb{R}\times [0,\infty)}e^{\mathrm{i}(u-w)x-uwy}\mu_{n_k}(dx,dy)\\
    &= \lim_{k\to\infty}\varphi_{n_k}(u)\varphi_{n_k}(-w) \\
    &= \varphi(u)\varphi(-w).
    \end{aligned}
    \]
    Therefore, $X\in \mathcal{H}_{BB}$ and since the Fourier-Laplace transform uniquely determines $\mu$ (see the latter part of the proof of Proposition~\ref{prop:Thorin_coupling}), the whole sequence $\mu_n$ converges to $\mu$.
\end{proof}

\begin{theorem}\label{thm:HCM_multiplication_closure}
    Let $X\in \mathcal{H}_{BB}$ and let $Y\in \mathcal{H}_{CM}$ be independent and strictly positive. Then $YX\in \mathcal{H}_{BB}$ and $X/Y\in \mathcal{H}_{BB}$.  
\end{theorem}

The Bochner-Bernstein measure $\mu$ and the law of the associated coupling $(Z,Y)$ have mostly been treated abstractly so far. The Bochner-Bernstein measure can generally only be computed in terms of the Thorin measure of $X$, but below we characterize some properties of $\mu$. Proposition~\ref{prop:coupling_type_G} shows that the second coordinate of the coupling $(Z,Y)$ is itself a generalized gamma convolution with marginal distribution that carries information on the symmetric part of the law of $X$. In particular, the symmetrization of $X\in\mathcal{T}(\mathbb{R})$ has a type-G representation with variance $2Y$ and the Thorin measure of $Y$ is the image under $x\mapsto x^2$ of the sum of the positive and negative Thorin measures of $X$, thereby generalizing Example~\ref{ex:laplace_coupling}. Proposition~\ref{prop:thorin_cumulants} states the joint cumulants of the coupling $(Z,Y)$ in terms of the cumulants of $X$.

\begin{proposition}\label{prop:coupling_type_G}
    Let $X\in\mathcal{T}(\mathbb{R})$ with characteristic function $\varphi$ of the form \eqref{eq:thorin_class_def} and Thorin measure $U$, and let $X'$ be an independent copy of $X$. Let $(Z,Y)$ denote the coupling of Proposition~\ref{prop:Thorin_coupling} of $X$, which exists due to Theorem~\ref{thm:Thorin_HFL_equivalent}. Then the following holds
    \begin{enumerate}
        \item $\mathbb{E}\left[e^{-qY}\right]=\lvert \varphi(\sqrt{q})\rvert^2$ for all $q\geq0$. In particular, the law of $Y$ is determined by $\lvert\varphi\rvert$ and is the same for every coupling $(Z,Y)$.
        \item $X-X'\overset{d}{=}\sqrt{2Y}N$, where $N\sim N(0,1)$ is independent of $Y$.
        \item $Y\in\mathcal{T}(\mathbb{R}_+)$ with Thorin measure $U_Y$ given by the image of $U$ under $x\mapsto x^2$ such that
        \[
        \int_0^\infty g(x)U_Y(dx)=\int_0^\infty g(x^2)\left( U_{+}+U_{-} \right)(dx).
        \]
        Additionally, $Y\geq \sigma^2$ almost surely.
    \end{enumerate}
\end{proposition}
\begin{proof}
    The first claim follows directly from \eqref{eq:mu_integral_bound}. For $q=0$, both sides of \eqref{eq:mu_integral_bound} are $1$ and since finite measures on $[0,\infty)$ are determined by their Laplace transforms on $(0,\infty)$, the law of $Y$ is determined by $q\mapsto \lvert \varphi(\sqrt{q})\rvert^2$, hence by $\lvert\varphi \rvert$ and does not depend on the coupling $(Z,Y)$. The second claim follows by setting $u=w$ in item three of Proposition~\ref{prop:Thorin_coupling}. We now show the third claim. Formally, setting $u=\sqrt{w}$,
    \[
    \begin{aligned}
    \mathbb{E}\left[ e^{-wY} \right] &= \exp\left( -\sigma^2w-\int_{\mathbb{R}\setminus \lbrace 0\rbrace}\log\left( 1+\frac{w}{x^2}\right)U(dx) \right) \\
    &= \exp\left( -\sigma^2 w-\int_0^\infty \log\left( 1+\frac{w}{y^2}\right)(U_++U_-)(dy) \right).
    \end{aligned}
    \]
    It follows that, if $\log(1+w/x^2)$ is integrable with respect to $U$, then $U_Y$ is a Thorin measure of a GGC, since
    \[
    \begin{gathered}
    \int_{(0,1]}\lvert \log(y)\rvert U_Y(dy)=\int_{\lbrace 0<\lvert x \rvert\leq 1\rbrace}\lvert \log(x^2)\rvert U(dx)<\infty, \\
    \int_{(1,\infty)}y^{-1}U_Y(dy)= \int_{\lbrace \lvert x\rvert >1\rbrace}x^{-2}U(dx)\leq 2\int_{\lbrace\lvert x\rvert >1\rbrace}\frac{1}{1+x^2}U(dx)<\infty.
    \end{gathered}
    \]
    $Y$ is thus a GGC with $a=\sigma^2$ and Thorin measure $U_Y$. The pair $(a,U_Y)$ uniquely defines the GGC by Theorem~3.1.1 in \textcite{Bondesson92} and since $U_Y$ is positive, it follows that $Y\geq a=\sigma^2$. We now argue that $\log(1+w/x^2)$ is integrable with respect to $U$, which completes the proof. Lemma~\ref{lem:U_integrability} below yields that $\log(x/(x-\mathrm{i}u))-\mathrm{i}ux/(1+x^2)$ is absolutely integrable with respect to $U$ and hence we may take real parts under the integral sign to obtain that $\log (1+u^2/x^2)$ is integrable with respect to $U$. Setting $u=\sqrt{w}$ completes the argument.
\end{proof}
\begin{proposition}\label{prop:thorin_cumulants}
    Let $X\in\mathcal{T}(\mathbb{R})$ with $\mathbb{E}[\lvert X\rvert^K]<\infty$ for some $K\geq 1$ and with cumulants $\kappa_k(X)$ for $k=1,\ldots ,K$. Let $(Z,Y)$ denote the coupling of Proposition~\ref{prop:Thorin_coupling} of $X$, which exists due to Theorem~\ref{thm:Thorin_HFL_equivalent}. Then the mixed moments $\mathbb{E}[\lvert Z\rvert^i Y^j]$ are finite with $i+2j\leq K$ and the joint cumulants of $(Z,Y)$ are given by
    \[
    \kappa_{i,j}(Z,Y) = \frac{(i+j-1)!}{(i+2j-1)!}\kappa_{i+2j}(X), \quad 1\leq i+2j\leq K.
    \]
\end{proposition}
\begin{proof}
    We first argue that the mixed moments $\mathbb{E}[\lvert Z\rvert^iY^j]$ are finite for $i+2j\leq K$. For $i,j\geq 1$, set $n=i+2j$ and note that Hölder's inequality gives
    \[
    \mathbb{E}\left[ \lvert Z\rvert^i Y^j \right] \leq \left( \mathbb{E}\left[ \lvert Z\rvert^n \right] \right)^{\frac{i}{n}}\left( \mathbb{E}\left[ Y^\frac{n}{2} \right]\right)^{\frac{2j}{n}} = \left( \mathbb{E}[\lvert X\rvert^n] \right)^{\frac{i}{n}}\left(\mathbb{E}\left[ Y^{\frac{n}{2}} \right] \right)^{\frac{2j}{n}},
    \]
    where the last equality follows from $Z\overset{d}{=}X$ due to Proposition~\ref{prop:Thorin_coupling}. For an independent copy $X'$ of $X$, Proposition~\ref{prop:coupling_type_G} yields that $X-X'\overset{d}{=}\sqrt{2Y}N$ with $N$ standard Gaussian and independent of $Y$, so $\mathbb{E}\left[\lvert X-X'\rvert\right]=2^{n/2}\mathbb{E}[Y^{n/2}]\mathbb{E}[\lvert N \rvert^n]$. Thus, $\mathbb{E}[Y^{n/2}]$ is finite exactly when $\mathbb{E}[\lvert X-X'\rvert^n]<\infty$, equivalently when $\mathbb{E}[\lvert X\rvert^n]<\infty$. The cases $i=0$ and $j=0$ yield the two factors on their own. Hence we find that finiteness of $\mathbb{E}[\lvert X\rvert^K]$ yields finiteness of $\mathbb{E}[\lvert Z\rvert^i Y^j]$ whenever $i+2j\leq K$ and the case $i=K,j=0$ shows that the assumption is also necessary. The same bound with $(i,j)$ replaced by $(n-l,l/2)$ shows that $\mathbb{E}\left[ \lvert Z\rvert^{n-l}Y^{l/2} \right]<\infty$ for all $0\leq l\leq n\leq K$. 

    Now fix $u\in\mathbb{R}$ and $w>0$ and let $N$ be a standard Gaussian random variable independent of $Z,Y$ and define
    $
    W = uZ + \sqrt{2w Y}N.
    $
    Conditionally on $(Z,Y)$, $W$ is then Gaussian with mean $uZ$ and variance $2wY$ such that $\mathbb{E}[e^{\mathrm{i}hW}]=\mathbb{E}[e^{\mathrm{i}huZ-h^2wY}]$. Set $S_\pm = (\pm u+\sqrt{u^2+4w})/2$ such that $S_+-S_- = u$ and $S_+S_- = w$. For $h>0$, it follows from \eqref{eq:HBB_joint_FL_characterization} with $p=hu$ and $q=h^2w$ that 
    \[
    \mathbb{E}\left[ e^{\mathrm{i}hW} \right] = \varphi(hS_+)\overline{\varphi(hS_-)}.
    \]
    $S_+$ and $(-S_-)$ are the roots of the polynomial $x^2-ux-w=0$ and Newton's identities therefore give $(S_+)^n+(-S_-)^n=P_n(u,w)$ where
    \[
    P_n(u,w)=\sum_{k\leq n/2}\frac{n}{n-k}\binom{n-k}{k}u^{n-2k}w^k.
    \]
    Since $\mathbb{E}[\lvert X\rvert^K]<\infty$, $\varphi_X$ is $K$ times continuously differentiable and the distinguished logarithm $\psi(z)=\log(\varphi_X(z))$ satisfies
    \[
    \psi (z) =  \sum_{n\leq K}\kappa_n(X)\frac{(\mathrm{i}z)^n}{n!}+o\left(\lvert z\rvert^K\right), \quad z\to 0.
    \]
    Conjugate symmetry of $\varphi_X$ yields $\overline{\psi(z)}=\psi(-z)$ for real $z$, and since $(\mathrm{i}hS_+)^n+(-\mathrm{i}hS_-)^n=(\mathrm{i}h)^n P_n(u,w)$, we thus find that for $h\to 0$,
    \[
    \log\left( \mathbb{E}\left[ e^{\mathrm{i}hW}\right] \right) = \psi(hS_+)+\psi(-hS_-) = \sum_{n\leq K}\frac{\kappa_n(X)(\mathrm{i}h)^nP_n(u,w)}{n!}+o\left( h^K\right).
    \]
    The $K$'th absolute moment of $W$ is also finite, since 
    \[
    \mathbb{E}[\lvert W\rvert^K]\leq 2^{K-1}\left(\lvert u\rvert^K \mathbb{E}[\lvert X\rvert^K]+(2w)^{K/2}\mathbb{E}[Y^{K/2}]\mathbb{E}[\lvert N\rvert^K] \right)<\infty,
    \]
    and thus it also holds that $\log\left(\mathbb{E}[e^{\mathrm{i}hW}]\right)=\sum_{n\leq K}\kappa_n(W)\tfrac{(\mathrm{i}h)^n}{n!}+o\left(h^K\right)$. Matching coefficients in the two expansions, it must therefore hold that $\kappa_n(W)=\kappa_n(X)P_n(u,w)$. Expanding $(uZ+\sqrt{2wY}N)^n$ and using that $N$ is independent of $(Z,Y)$ with $\mathbb{E}[N^l]=0$ for odd $l$ and $\mathbb{E}[N^{2k}]=(2k)!/(2^kk!)$ yields
    \[
    \mathbb{E}\left[W^n\right] = \sum_{i+2j=n}\frac{n!}{i!j!}u^i w^j\mathbb{E}\left[ Z^iY^j \right], \quad n \leq K.
    \]
    Therefore,
    \[
    \kappa_n(W) = \sum_{i+2j=n}\frac{n!}{i!j!}u^i w^j \kappa_{i,j}(Z,Y).
    \]
    So for each $n\leq K$, the two expressions for $\kappa_n(W)$ are polynomials in $(u,w)$ that agree for all $u\in\mathbb{R}$ and $w> 0$ and hence their coefficients must match. We therefore find that
    \[
    \kappa_{i,j}(Z,Y)=\kappa_n(X) \frac{i! j!}{n!}\frac{n}{n-j}\frac{(n-j)!}{j!i!} = \kappa_n(X) \frac{(n-j-1)!}{(n-1)!},
    \]
    which is the stated formula.
\end{proof}

\section{Applications to laws of exponential functionals of Lévy processes}\label{sec:applications}
Proposition~\ref{prop:recursion} below is the main result of the second part of the paper and it implies membership of $\mathcal{T}(\mathbb{R})$ for a large class of skewed and heavy-tailed distributions. The laws of Proposition~\ref{prop:recursion} correspond, e.g., to those of insurance risk models with stochastic return on investments similar to those of \textcite{Paulsen1993,Norberg1999,Paulsen2008}. Related results in the symmetric setting are given in \textcite{BehmeMaejimaMatsuiSakuma2012,Bondesson2015} and some considerations on the asymmetric and light-tailed case are also given in \textcite{Torricelli2026}. The ingredient that permits us to treat the more general skewed and possibly heavy-tailed distributions is exactly Corollary~\ref{cor:HCM_closure_Thorin}. Having established this result, we use it in Section~\ref{sec:ADY} to show that the class of Markov processes studied in \textcite{AliliDufresneYor1997} in the context of a generalized Bougerol identity has, at all times, its law in $\mathcal{T}(\mathbb{R})$. We then show in Example~\ref{ex:JSU} and Example~\ref{ex:Pearson} that this implies that Johnson's S\textsubscript{U} distribution \citep[][]{Johnson1949} and the Pearson type IV distributions are in $\mathcal{T}(\mathbb{R})$ for all admissible parametrizations. For both of these distribution classes, it is relatively simple to show via classical methods that the symmetric members are in $\mathcal{T}(\mathbb{R})$ due to their Gaussian mixture representations; for the Pearson type IV distribution this is simply the student's t-distribution and for Johnson's S\textsubscript{U} distribution the proof is outlined in Remark~\ref{rem:JSU_symmetric}. Section~\ref{sec:mixture} shows that a Gaussian mean-variance mixture representation fails for the asymmetric processes in the class of \textcite{AliliDufresneYor1997} with $b\neq 0$, and that the subclass with $\nu = 0$ consists instead mixtures of skewed generalized hyperbolic laws. 

\begin{proposition}\label{prop:recursion}
    Let $(B_t)_{t\geq 0}$ be a standard Brownian motion and let $a>0$, $m\in \mathbb{R}$, and $U_t = aB_t + mt$. Let $(L_t)_{t\geq 0}$ be a Lévy process that is independent of $B$ and with $L_1\in\mathcal{T}(\mathbb{R})$ and define for $z\in\mathbb{R}$ the processes
    \[
    Z_t = e^{U_t}\left( z+\int_0^t e^{-U_s}dL_s \right), \quad V_t = \int_0^t e^{U_s}dL_s.
    \]
    Then the following holds.
    \begin{enumerate}
        \item $Z_t,V_t\in \mathcal{T}(\mathbb{R})$ for every $t>0$ and $z\in\mathbb{R}$.
        \item If $V_t$ converges almost surely as $t\to \infty$ to a random variable $V_\infty$, then $V_\infty\in\mathcal{T}(\mathbb{R})$. 
    \end{enumerate}
\end{proposition}
\begin{proof}
    Fix $t>0$ and $n\in\mathbb{N}$ and set $\Delta=t/n$ and $t_k=k\Delta$ for $k=0,1,\ldots ,n$. Define 
    \[
    m_k = e^{U_{t_k}-U_{t_{k-1}}} = \exp\left( a(B_{t_k}-B_{t_{k-1}})+m\Delta \right), \quad \Delta L_k = L_{t_k}-L_{t_{k-1}}. 
    \]
    Set also $Z_0=z$ and $Z_k= m_k(Z_{k-1}+\Delta L_{k})$. A simple induction in $k$ shows that
    \begin{equation}\label{eq:Z_riemann_sum}
    Z_k = e^{U_{t_k}}\left( z +\sum_{j=1}^{k}e^{-U_{t_{j-1}}} \Delta L_j\right).
    \end{equation}
    Likewise, define $W_{n+1}=0$ and $W_k=\Delta L_k+m_k W_{k+1}$ for $k=n,(n-1),\ldots ,1$. Then a similar downwards induction shows that 
    \begin{equation}\label{eq:W_riemann_sum}
    W_k = \sum_{j=k}^{n}e^{U_{t_{j-1}}-U_{t_{k-1}}}\Delta L_j,\quad W_1 =\sum_{j=1}^{n}e^{U_{t_{j-1}}}\Delta L_j,
    \end{equation}
    such that $Z_n$ is the left-point Riemann sum approximation of $Z_t$ and $W_1$ is the left-point Riemann sum approximation of $V_t$. Let $\mathcal{G}_k=\sigma(B_s,L_s:s\leq t_k)$ be the $\sigma$-algebra generated by $B,L$ up to time $t_k$. Then $Z_{k-1}$ is $\mathcal{G}_{k-1}$-measurable and $m_k$ is a function of $B_{t_k}-B_{t_{k-1}}$ and $\Delta L_k$ is a function of $L_{t_k}-L_{t_{k-1}}$. As $B,L$ have independent increments and are mutually independent, the three random variables $Z_{k-1},m_k,\Delta L_k$ are all mutually independent. Similarly, $W_{k+1}$ is a function of $(B_s-B_{t_k},L_{s}-L_{t_k})_{s\geq t_k}$ so $W_{k+1},\Delta L_k, m_k$ are mutually independent. Note also that $L_\Delta\in\mathcal{T}(\mathbb{R})$ because $L_1\in\mathcal{T}(\mathbb{R})$ and that $\log (m_k)\sim N(m\Delta ,a^2\Delta)$ with $a^2\Delta >0$, meaning that $m_k\in\mathcal{H}_{CM}$ due to Example~5.2.1 of \textcite{Bondesson92}. Since degenerate laws are also in $\mathcal{T}(\mathbb{R})$ we have that $Z_0\in\mathcal{T}(\mathbb{R})$ and thus Corollary~\ref{cor:HCM_closure_Thorin} combined with closure of $\mathcal{T}(\mathbb{R})$ under convolution of independent random variables and an induction argument shows that $Z_k=m_k(Z_{k-1}+\Delta L_k)\in\mathcal{T}(\mathbb{R})$ for all $k=1,\ldots ,n$. A similar argument shows that $W_k\in\mathcal{T}(\mathbb{R})$ for all $k=1,\ldots,n$. 

    Letting $n\to\infty$, $L$ is a Lévy process in the filtration $\mathcal{G}_t=\sigma(L_s,B_s:s\leq t)$ and the maps $s\mapsto e^{-U_s}$ and $s\mapsto e^{U_s}$ are continuous and $\mathcal{G}_t$-adapted processes. The probability limits of the Riemann sums \eqref{eq:Z_riemann_sum} and \eqref{eq:W_riemann_sum} are
    \[
    \sum_{j=1}^{n}e^{-U_{t_{j-1}}}\Delta L_j \xrightarrow[n\to\infty]{\mathbb{P}}\int_0^t e^{-U_s}dL_s, \quad W_1 \xrightarrow[n\to\infty]{\mathbb{P}}V_t.
    \]
    Multiplying the first sum by $e^{U_t}$ and adding $e^{U_t}z$ yields that $Z_n\to Z_t$ in probability. Since convergence in probability implies convergence in law and $\mathcal{T}(\mathbb{R})$ is closed under weak limits, we conclude that $Z_t,V_t\in\mathcal{T}(\mathbb{R})$. 

    For the second statement, note that if $V_t\to V_\infty$ almost surely, then $V_t\to V_\infty$ also in law and so closure under weak limits of $\mathcal{T}(\mathbb{R})$ yields that $V_\infty\in\mathcal{T}(\mathbb{R})$. 
\end{proof}

\subsection{Recovering some known classes of distributions}\label{sec:ADY}
\begin{corollary}\label{cor:ADY_is_thorin}
    Let $\mu,\nu,b\in\mathbb{R}$ be constants, $(\beta_t)_{t\geq 0}$ a Brownian motion that is independent of $B$, and let $Y$ be the diffusion process with generator $\tfrac{1}{2}\tfrac{d^2}{dy^2}+(\mu\tanh(y)+\nu\operatorname{sech}(y))\tfrac{d}{dy}$ and starting in $Y_0=b$. Then $(Z_t)_{t\geq 0}\overset{d}{=}(\sinh(Y_t))_{t\geq 0}$ where 
    \begin{equation}\label{eq:ADY_Z}
    Z_t = e^{B_t+\mu t}\left( \sinh(b)+\int_0^t e^{-(B_s+\mu s)}(d\beta_s+\nu ds)\right).
    \end{equation}
    In particular, $\sinh(Y_t)\in\mathcal{T}(\mathbb{R})$ for all $t>0$.
\end{corollary}
\begin{proof}
    The fact that $Z_t\overset{d}{=}\sinh(Y_t)$ can be deduced from Proposition~2.1 of \textcite{Vakeroudis2012} \citep[attributed to][]{AliliDufresneYor1997} applied with their $\mu,\nu$ and starting point $x=\sinh(b)$. Setting $z=\sinh(b)$, $a=1$, $m=\mu$, $L_t=\beta_t+\nu t$ in Proposition~\ref{prop:recursion} shows that $Z_t\in\mathcal{T}(\mathbb{R})$ and hence $\sinh(Y_t)\in\mathcal{T}(\mathbb{R})$ for all $t>0$.  
\end{proof}

\begin{corollary}\label{cor:A_t_is_GGC}
    Let $\alpha\in \mathbb{R}$ and define the process
    \(
    A_t^{(\alpha)} = \int_0^t e^{2(B_s+\alpha s)}ds.
    \)
    Then $A_t^{(\alpha)}\in \mathcal{T}(\mathbb{R}_+)$.
\end{corollary}
\begin{proof}
    It follows from Proposition~\ref{prop:recursion} by setting $a=2$, $m=2\alpha$, and $L_t=t$ that $A_t^{(\alpha)}\in\mathcal{T}(\mathbb{R})$. Since $A_t^{(\alpha)}$ is almost surely positive, it then follows that $A_t^{(\alpha)}\in\mathcal{T}(\mathbb{R}_+)$. 
\end{proof}

\begin{remark}\label{rem:A_t_is_GGC}
    All quantities in Corollary~\ref{cor:A_t_is_GGC} are positive. This result can thus be proved using only the previously established multiplicative closure of $\mathcal{T}(\mathbb{R}_+)$ by re-doing the proof of Proposition~\ref{prop:recursion} for $\mathcal{T}(\mathbb{R}_+)$. As such, this result does not rely on Corollary~\ref{cor:HCM_closure_Thorin} or the hyperbolic Bochner-Bernstein property of the laws in $\mathcal{T}(\mathbb{R})$. The law of $A_t^{(\alpha)}$ has been studied by several authors; \textcite{Dufresne1990} establishes that the perpetuity $A_\infty^{(\alpha)}$ has an inverse gamma distribution when $\alpha <0$ and \textcite{JedidiVakeroudis2018} show that $1/A_t^{(\alpha)}\in\mathcal{T}(\mathbb{R}_+)$ for $\alpha=0$. To the best of our knowledge, it has not previously been explicitly stated in the literature that $A_t^{(\alpha)}\in\mathcal{T}(\mathbb{R}_+)$ for fixed $t>0$ and general $\alpha$.
\end{remark}
\begin{remark}\label{rem:JSU_symmetric}
    Corollary~\ref{cor:ADY_is_thorin} shows that the law of the process $\sinh(b+B_t)$ is in $\mathcal{T}(\mathbb{R})$ at all times (set $\mu=\nu=0$). The symmetric result with $b=0$ follows directly from Bougerol's identity \citep[][]{Bougerol1983}, which shows that $\sinh(B_t)\overset{d}{=}\sqrt{A_t^{(0)}}N$ with $A_t^{(0)}$ as in Corollary~\ref{cor:A_t_is_GGC}. Indeed, Chapter~7 of \textcite{Bondesson92} establishes that symmetric Thorin laws are Gaussian mixtures with GGC mixing variables and that the converse also holds, so $\sinh(B_t)\in\mathcal{T}(\mathbb{R})$. By Remark~\ref{rem:A_t_is_GGC}, the fact that $A_t^{(0)}\in\mathcal{T}(\mathbb{R}_+)$ can be proven from just the well-known results of \textcite{Bondesson92}, so to show $\sinh(B_t)\in\mathcal{T}(\mathbb{R})$ we do not require the HBB property. The asymmetric case with $b\neq 0$, however, is not susceptible to this approach, but Corollary~\ref{cor:HCM_closure_Thorin} (and thus the HBB property) provides a feasible route.
\end{remark}

\begin{example}[Johnson's S\textsubscript{U} distribution]\label{ex:JSU}
    Johnson's S\textsubscript{U} distribution was introduced in \textcite{Johnson1949} and is typically defined as the law of the random variable
    \[
    Y = \xi + \delta \sinh\left( \frac{\lambda+N}{\sigma} \right) = \xi+\frac{\delta\left( e^{(\lambda+N)/\sigma} -e^{-(\lambda+N)/\sigma} \right)}{2},
    \]
    where $\xi,\lambda\in\mathbb{R}$, $\delta,\sigma>0$ and $N$ is a standard Gaussian random variable. This is frequently used as a two-sided version of the lognormal distribution, with which it shares many properties such as moments of all orders and moment indeterminacy. The distribution is frequently employed as a flexible and skewed distribution for heavy-tailed real-valued data and has been utilized within diverse areas such as economics, risk modelling, and sea clutter modelling; see for example \textcite{JSU_commodities,JSU_futures,JSU_sea_clutter}. It follows by setting $\mu=\nu=0$ in Corollary~\ref{cor:ADY_is_thorin} and closure of $\mathcal{T}(\mathbb{R})$ under affine maps that Johnson's S\textsubscript{U} distribution is in $\mathcal{T}(\mathbb{R})$ for any choice of parameters.
\end{example}
\begin{example}[Pearson type IV distribution]\label{ex:Pearson}
    In the notation of Proposition~\ref{prop:recursion}, let $\mu<0$, $\nu,b\in\mathbb{R}$, $U_t = B_t+\mu t$, $L_t=\beta_t+\nu t$ with $(\beta_t)_{t\geq 0}$ a standard Brownian motion, independent of $B$. Then \eqref{eq:ADY_Z} reads
    \[
    Z_t = e^{U_t}\sinh (b) +\int_0^t e^{U_t-U_s}(d\beta_s + \nu ds) \overset{d}{=}\sinh(Y_t).
    \]
    For fixed $t$, the reversed pair $(U_t-U_{t-s},L_t-L_{t-s})_{0\leq s\leq t}$ has the law of $(U_s,L_s)_{0\leq s\leq t}$, since $(U,L)$ has stationary independent increments and continuous paths. Hence $\int_0^t e^{U_t-U_s}dL_s \overset{d}{=}V_t$. Theorem~2.1 of \textcite{KondoMaejimaSato2006} yields that $V_t\to V_\infty=\int_0^\infty e^{U_s}dL_s$ almost surely and since $e^{U_t}\sinh(b)\to 0$ almost surely, we have $Z_t\to V_\infty$ in law for any starting value $b$. The density of $V_\infty$ is given in Section~6.1 of \textcite{BertoinYor2005} as $f(x)=C(1+x^2)^{\mu-1/2}e^{2\nu \arctan(x)}$. By the convergence just shown, this is also the stationary law of $\sinh(Y_t)$. This is a Pearson type IV distribution with location parameter $0$, scale parameter $1$, and shape parameters $(1/2-\mu)>1/2$ and $2\nu$. Self-decomposability of $V_\infty$ then follows from Theorem~3.1 of \textcite{KondoMaejimaSato2006}, but Corollary~\ref{cor:ADY_is_thorin} yields that $\sinh(Y_t)\in\mathcal{T}(\mathbb{R})$ and closure under weak limits therefore yields the stronger statement that $V_\infty\in\mathcal{T}(\mathbb{R})$. Closure of $\mathcal{T}(\mathbb{R})$ under affine maps then shows that any Pearson type IV distribution is in $\mathcal{T}(\mathbb{R})$ (equivalently, apply the second statement of Proposition~\ref{prop:recursion}). The case $\nu=0$ recovers the well-known fact that the student's t-distribution is in $\mathcal{T}(\mathbb{R})$.
\end{example}

\subsection{Mixture representations}\label{sec:mixture}
Lemma~\ref{lem:gaussian_mixture} provides general criteria for failure of a Gaussian mixture representation, whenever the distribution admits a continuous and strictly positive density. The criteria can loosely be stated as the density having a sub-linear logarithmic ratio and no center of symmetry. This is particularly relevant for the skewed laws of $\sinh(Y_t)$ arising when $b\neq 0$ in Corollary~\ref{cor:ADY_is_thorin}. Corollary~\ref{cor:gaussian_mixture} below uses Lemma~\ref{lem:gaussian_mixture} to show that the asymmetric distributions of Corollary~\ref{cor:ADY_is_thorin} with $b\neq 0$ are not Gaussian mean-variance mixtures and Corollary~\ref{cor:pearson_iv_mixture} with $\nu\neq 0$ shows this for the Pearson type IV distributions. Many tractable infinitely divisible distributions such as the generalized hyperbolic distributions of \textcite{Barndorff-NielsenHalgreen1977} are Gaussian mixtures. Proposition~\ref{prop:hyperbolic_mixture} then establishes that the law of $Z_t$ with $b\neq 0$ and $\nu=0$ in Corollary~\ref{cor:ADY_is_thorin} is in fact a mixture of such generalized hyperbolic distributions with a rather complicated mixing law.
\begin{lemma}\label{lem:gaussian_mixture}
    Let $X$ be a random variable with a continuous and strictly positive density $f$ on $\mathbb{R}$. Suppose that the following holds.
    \begin{enumerate}
        \item For every $h\in\mathbb{R}$, $\log(f(h+x))-\log(f(h-x))=o(x)$ as $x\to\infty$.
        \item There is no $h\in\mathbb{R}$ such that $f(h+x)=f(h-x)$ for all $x\in\mathbb{R}$. 
    \end{enumerate}
    Then $X$ is not a normal mean-variance mixture, i.e., there are no parameters $\ell,m\in\mathbb{R}$ and no non-negative random variable $V$ such that
    \begin{equation}\label{eq:gaussian_mixture}
    X\overset{d}{=}\ell + mV+\sqrt{V}N,
    \end{equation}
    where $N$ is a standard Gaussian random variable, independent of $V$.
\end{lemma}
\begin{proof}
    Suppose that a representation of the form \eqref{eq:gaussian_mixture} exists. Then, since $X$ has a density, it holds that $\mathbb{P}(V=0)=0$ and letting $\eta$ denote the law of $V$ on $(0,\infty)$ it holds that
    \[
    f(\ell +x)=\int_0^\infty \frac{1}{\sqrt{2\pi v}}e^{-\frac{(x-mv)^2}{2v}}\eta (dv) = e^{mx}g(x),
    \]
    where 
    \[
    g(x) = \int_0^\infty \frac{1}{\sqrt{2\pi v}}e^{-\frac{x^2}{2v}-\frac{m^2v}{2}}\eta (dv).
    \]
    Note that $g$ is even and positive, so $\log(f(\ell+x))-\log(f(\ell-x))=2mx$ for all $x$. If $m\neq 0$, setting $h=\ell$ contradicts the fact that $\log(f(h+x))-\log(f(h-x))=o(x)$ as $x\to\infty$ and if $m=0$, it contradicts $f(h+x)=f(h-x)$. We therefore conclude that $X$ cannot admit a Gaussian mixture representation of the form \eqref{eq:gaussian_mixture}.  
\end{proof}
\begin{corollary}\label{cor:gaussian_mixture}
    Let the setting be as in Corollary~\ref{cor:ADY_is_thorin} and suppose that $b\neq 0$. Then the law of $Z_t$ does not admit a Gaussian mixture representation of the form \eqref{eq:gaussian_mixture}. 
\end{corollary}
\begin{proof}
    By Corollary~\ref{cor:ADY_is_thorin}, it holds that $Z_t\overset{d}{=}\sinh(Y_t)$, so it suffices to study the law of $\sinh(Y_t)$. Set $g(y)=\mu \tanh(y)+\nu \operatorname{sech}(y)$ which is bounded and with bounded derivative $g'(y)=\mu \operatorname{sech}(y)^2-\nu \operatorname{sech}(y)\tanh(y)$. Let $F(y)=\mu \log(\cosh(y))+\nu \arctan(\sinh(y))$ and denote by $\mathbb{P}_b$ the Wiener measure of a Brownian motion $X$ started at $b$, and set 
    \[
    M_t = \exp\left( \int_0^t g(X_s)dX_s -\frac{1}{2}\int_0^t g(X_s)^2ds \right).
    \]
    Since $g$ is bounded, $M_t$ is a martingale and Girsanov's theorem yields that under $M_t\cdot \mathbb{P}_b$, the process $X$ solves the SDE $dX_t=dW_t+g(X_t)dt$ with $W_t$ a $M_t\cdot \mathbb{P}_b$-Brownian motion. It therefore holds that $\mathbb{E}[\phi (Y_t)]=\mathbb{E}^{\mathbb{P}_b}[\phi (X_t)M_t]$. Applying Itô's formula, we find that
    \[
    \int_0^t g(X_s)dX_s = F(X_t)-F(b)-\frac{1}{2}\int_0^t g'(X_s)ds,
    \]
    such that $M_t=\exp\left( F(X_t)-F(b)-\frac{1}{2}\int_0^t(g'(X_s)+g(X_s)^2)ds \right)$. Conditioning on $X_t=x$, we find that the density of $Y_t$ is
    \[
    f_{Y_t}(x)=p_t(x-b)e^{F(x)-F(b)}\mathbb{E}^{\mathbb{P}_b}\left[ e^{-\frac{1}{2}\int_0^t (g'(X_s)+g(X_s)^2)ds}\mid X_t =x \right].
    \]
    The expectation on the right hand side is in $[e^{-Kt/2},e^{Kt/2}]$, uniformly in $x$ with $K=\sup_{y}\lvert g'(y)+g(y)^2\rvert<\infty$. Note then that $\arctan(\sinh(x))$ is bounded in $x$ and that
    \[
    \log (\cosh(x))=\lvert x\rvert - \log (2) + \mathcal{O}(e^{-2\lvert x\rvert}), \quad \log(f_{Y_t}(x))=-\frac{(x-b)^2}{(2t)}+\mu\lvert x\rvert + \mathcal{O}(1),
    \]
    both uniformly in $x$ as $x\to\infty$. For $X=\sinh(Y_t)$, the density is $f_X(z)=f_{Y_t}(\operatorname{asinh}(z))/\sqrt{1+z^2}$ with $\operatorname{asinh}(z)=\operatorname{sign}(z)\log(2\lvert z\rvert) + \mathcal{O}(z^{-2})$, which is continuous and strictly positive. Thus we find that
    \begin{equation}\label{eq:log_f_X}
    \log(f_X(h+x))-\log(f_X(h-x)) = \frac{2b}{t}\log(2x)+\mathcal{O}(1),
    \end{equation}
    as $x\to\infty$, which is $o(x)$ such that the first condition of Lemma~\ref{lem:gaussian_mixture} holds. The right hand side of \eqref{eq:log_f_X} is also unbounded for $b\neq 0$, so $f_X$ cannot be symmetric about any $h$. Indeed, if $f_X$ is symmetric about some $h$, this would make the left hand side vanish and thus the second condition of Lemma~\ref{lem:gaussian_mixture} holds. We conclude that the law of $\sinh(Y_t)$ cannot be represented as a Gaussian mixture.
\end{proof}
\begin{corollary}\label{cor:pearson_iv_mixture}
    Let $X$ be a random variable with a Pearson type IV density of the form $f(x)=C(1+x^2)^{\mu-1/2}e^{2\nu\arctan(x)}$ with $\nu\neq 0$. Then $X$ does not admit a Gaussian mixture representation of the form \eqref{eq:gaussian_mixture}.
\end{corollary}
\begin{proof}
    Note that 
    \[
    \log(f(h+x))-\log (f(h-x)) = \left(\mu -\frac{1}{2}\right)\log\left( \frac{1+(h+x)^2}{1+(h-x)^2} \right) + 2\nu(\arctan(h+x)-\arctan(h-x)).
    \]
    The right hand side tends to $2\pi\nu$ as $x\to\infty$, which is $o(x)$ and non-zero when $\nu\neq 0$. Applying Lemma~\ref{lem:gaussian_mixture} concludes the proof.
\end{proof}
\begin{proposition}\label{prop:hyperbolic_mixture}
    Let the setting be as in Corollary~\ref{cor:ADY_is_thorin} and set $\nu=0$. Define $\mathcal{Y}_t = e^{U_t}\int_0^t e^{-2U_s}ds$ with $U_t=B_t+\mu t$. Then, conditionally on $\mathcal{Y}_t = y$, we have
    \begin{equation}\label{eq:Z_mixture}
    Z_t \overset{d}{=}\frac{\sinh(b)}{y}V+\sqrt{V}N, 
    \end{equation}
    where $V\sim GIG(\mu,1,y^{-2})$ where $GIG(x_1,x_2,x_3)$ denotes the generalized inverse Gaussian distribution with density proportional to $v^{x_1-1}\exp\left(-\frac{1}{2}(\frac{x_2}{v}+vx_3 )\right)$ and $N\sim N(0,1)$ is independent of $V$ and $B$.
\end{proposition}
\begin{proof}
    Fix $t>0$. Conditionally on $B$, the process $Z_t$ is Gaussian with mean $\sinh(b)e^{U_t}$ and variance $\widetilde{A}_t=\int_0^t e^{2(U_t-U_s)}ds$. Set $G=e^{U_t}$ so $\widetilde{A}_t=G\mathcal{Y}_t$ pathwise and, conditionally on $\mathcal{Y}_t=y$ and $G=g$, $Z_t\sim N(\sinh(b)g,yg)$. The reversed process $\widetilde{U}_s=U_t-U_{t-s}$ is again a Brownian motion with drift $\mu$ and $\widetilde{U}_t=U_{t}$, so $\widetilde{A}_t=\int_0^t e^{2\widetilde{U}_s}ds$ and 
    \[
    (G,\mathcal{Y}_t)\overset{d}{=}\left( e^{U_t}, e^{-U_t}A_t^{(\mu)} \right),
    \]
    with $A_t^{(\mu)}$ as defined in Corollary~\ref{cor:A_t_is_GGC}. The conditional law of $G$ given $\mathcal{Y}_t=y$ is that of $e^{U_t}$ given $e^{-U_t}A_t^{(\mu)}=y$, which is $GIG(\mu,1/y,1/y)$ due to Proposition~5.2 of \textcite{MatsumotoYor2005_2}. This yields \eqref{eq:Z_mixture}. 
\end{proof}

\section{Proofs of the main results}\label{sec:proofs}
This section gives the proofs of Theorem~\ref{thm:HCM_multiplication_closure} and Theorem~\ref{thm:Thorin_HFL_equivalent}, which relies on some technical lemmas. Lemma~\ref{lem:pick_theorem} is essentially Theorem~7.1.2 of \textcite{Bondesson92}, which also shows the converse direction. We state it here for completeness and in our notation to make the required assumptions explicit. 
\begin{lemma}\label{lem:U_integrability}
    Let $U$ be a positive measure on $\mathbb{R}\setminus \lbrace 0 \rbrace$ satisfying the integrability conditions \eqref{eq:Thorin_measure_integrability}. For real $r,u,x$ with $x\neq 0$, define
    \[
    \begin{gathered}
    f_r(x)=\frac{1}{x-\mathrm{i}r}-\frac{x}{1+x^2}=\frac{1+\mathrm{i}xr}{(x-\mathrm{i}r)(1+x^2)}, \\
    g_u(x)=\log\left( \frac{x}{x-\mathrm{i}u} \right) - \frac{\mathrm{i}ux}{1+x^2}=-\log \left( 1-\frac{\mathrm{i}u}{x}\right)-\frac{\mathrm{i}ux}{1+x^2},
    \end{gathered}
    \]
    where the principal logarithm is used. Then the following holds.
    \begin{enumerate}
        \item $\int_{\mathbb{R}\setminus \lbrace 0\rbrace} \lvert g_u(y)\rvert U(dy)<\infty$.
        \item $\int_0^{\lvert u\rvert}\int_{\mathbb{R}\setminus \lbrace 0\rbrace}\lvert f_r(y)\rvert U(dy)dr<\infty$.
        \item $U$ is finite on $\lbrace 0<\lvert x\rvert \leq 1\rbrace$. In particular, $\int \left( \mathbf{1}_{\lbrace \lvert y\rvert \geq 1\rbrace}(y)\tfrac{1}{y}-\tfrac{y}{1+y^2} \right) U(dy)$ is finite.
    \end{enumerate}
\end{lemma}
\begin{proof}
    On the segment from $0$ to $u$, the function $r\mapsto x-\mathrm{i}r$ does not vanish, so $\log (x-\mathrm{i}r)$ has a continuous branch there, with derivative $-\mathrm{i}/(x-\mathrm{i}r)$. Hence $\int_0^u(x-\mathrm{i}r)^{-1}dr=\mathrm{i}\ell(u)$ with $\ell(u)=\log(x-\mathrm{i}u)-\log(x)$ along that branch. It is straightforward to show that $\ell(u)=\log (1-\mathrm{i}u/x)$ where $\log$ denotes the principal logarithm. We thus find that
    \[
    \mathrm{i}\int_0^u f_r(x)dr = -\log\left( 1-\frac{\mathrm{iu}}{x} \right) - \frac{\mathrm{i}ux}{1+x^2}=g_u(x).
    \]
    Note that $f_{-r}(x)=\overline{f_r(x)}$, so $\lvert g_u(x)\rvert\leq \int_0^{\lvert u\rvert}\lvert f_{\pm r}(x)\rvert dr$. We now proceed to bound the right hand side. To this end, note that for $\lvert x\rvert \geq 1$, it holds that $\lvert 1+\mathrm{i}xr \rvert=\sqrt{1+x^2r^2}\leq \sqrt{1+x^2}\sqrt{1+r^2}$ and $\lvert x\rvert\sqrt{1+x^2}\geq (1+x^2)/\sqrt{2}$, so $\lvert f_r(x)\rvert \leq \sqrt{2}\sqrt{1+r^2}(1+x^2)^{-1}$. Note also that $g_u(x)=\mathrm{i}\int_0^u f_r(x)dr$ and set $D_u(x)=\int_0^{\lvert u\rvert}\lvert f_r(x)\rvert dr$. Then it follows that 
    \[
    D_u(x)\leq \sqrt{2}\lvert u \rvert\sqrt{1+u^2}(1+x^2)^{-1}\leq 2 \lvert u\rvert(1+\lvert u\rvert)(1+x^2)^{-1}.
    \]
    For $0<\lvert x\rvert \leq 1$ we have $\lvert f_r(x)\rvert \leq (x^2+r^2)^{-1/2}+\tfrac{1}{2}$, so 
    \begin{equation}\label{eq:D_bound_lemma}
    D_u(x) \leq \operatorname{asinh}(\lvert u\rvert/\lvert x\rvert) + \frac{\lvert u \rvert}{2} \leq \log\left( 1+\frac{2\lvert u\rvert}{\lvert x\rvert} \right) + \frac{\lvert u\rvert}{2}\leq \log(1+2\lvert u\rvert)+\lvert \log (\lvert x\rvert)\rvert+ \frac{\lvert u\rvert}{2}.
    \end{equation}
    The right hand side of \eqref{eq:D_bound_lemma} is at most $(2\log(1+2\lvert u\rvert)+\lvert u \rvert)(1+x^2)^{-1}+\lvert \log (\lvert x\rvert)\rvert$, and we thus obtain the bound
    \begin{equation}\label{eq:D_bound_2_lemma}
    \lvert g_u(x)\rvert \leq D_u(x) \leq \frac{C(u)}{1+x^2}+\mathbf{1}_{\lbrace \lvert x\rvert \leq 1\rbrace}\lvert \log (\lvert x\rvert)\rvert, 
    \end{equation}
    where $C(u)=2\log(1+2\lvert u\rvert)+\lvert u\rvert + 2\lvert u \rvert(1+\lvert u\rvert)$. The right hand side is integrable with respect to $U$ due to the integrability conditions \eqref{eq:Thorin_measure_integrability}. Hence we conclude that $g_u$ is absolutely integrable with respect to $U$. The second point follows from Tonelli's theorem applied to the non-negative function $(r,x)\mapsto \lvert f_r(x)\rvert$ on $[0,\lvert u \rvert]\times (\mathbb{R}\setminus \lbrace 0\rbrace)$. This gives $\int_0^{\lvert u\rvert}\int \lvert f_r(x)\rvert U(dx)dr = \int \lvert D_u(x)\rvert U(dx)$, which is finite due to \eqref{eq:D_bound_2_lemma}. Fubini's theorem therefore permits interchanging the integrals and the inner integral is $\mathrm{i}\int_0^u f_r(x)dr=g_u(x)$. For the third claim, note that for $\varepsilon<1$, it holds that $U(\lbrace \varepsilon\leq \lvert x\rvert \leq 1\rbrace)\leq 2\int_{\lbrace \varepsilon\leq \lvert x\rvert \leq 1\rbrace} (1+x^2)^{-1}U(dx)$ and $U(\lbrace 0<\lvert x\rvert <\varepsilon\rbrace)\leq (2\log (1/\varepsilon))^{-1}\int_{\lbrace 0< \lvert x\rvert \leq 1\rbrace}\lvert \log (x^2)\rvert U(dx)$. Both of these integrals are finite due to \eqref{eq:Thorin_measure_integrability} and the stated integral is finite since the integrand is bounded by $\tfrac{1}{2}$ for $\lvert x\rvert < 1$ and $\tfrac{1}{(1+x^2)}$ for $\lvert x\rvert \geq 1$. 
\end{proof}

\begin{lemma}\label{lem:pick_theorem}
    Let $\varphi$ be a characteristic function whose restriction on $(0,\infty)$ extends analytically to $\mathcal{H}=\lbrace u\in\mathbb{C}:\Re (u)>0\rbrace$ and without zeros in $\mathcal{H}$. Suppose that $\Re(\varphi'(u)/\varphi(u))\leq 0$ for all $u\in\mathcal{H}$. Then $\varphi$ is the Fourier transform of a law in $\mathcal{T}(\mathbb{R})$. 
\end{lemma}
\begin{proof}
    As $\mathcal{H}$ is simply connected and $\varphi$ has no zeros in $\mathcal{H}$, $\varphi$ admits a holomorphic logarithm $\psi=\log(\varphi)$ on $\mathcal{H}$ which is such that $\lim_{u\to 0}\psi(u)=0$. Define the function $P$ as
    $
    P(z)=-\mathrm{i}\psi'(-\mathrm{i}z),
    $
    for $\Im(z)>0$ and note that $\Im(z)>0$ means that $-\mathrm{i}z\in\mathcal{H}$. Then $P$ is holomorphic on the upper half-plane and $\Im (P(z))=-\Re (\psi'(-\mathrm{i}z))=-\Re(\varphi'(-\mathrm{i}z)/\varphi(-\mathrm{i}z))$. Thus $P$ is a Pick function in the sense of \textcite{Bondesson92} and it therefore has Nevanlinna representation
    \[
    P(z)=\beta +\sigma^2 z +\int_{\mathbb{R}}\left( \frac{1}{x-z}-\frac{x}{1+x^2} \right)U(dx),
    \]
    for some $\beta\in\mathbb{R}$, $\sigma^2\geq 0$ and a positive measure $U$ on $\mathbb{R}$ satisfying $\int_\mathbb{R}\frac{1}{1+x^2}U(dx)<\infty$. Using that $\psi'(u)=\mathrm{i}P(\mathrm{i}u)$, we find that, for $u>0$,
    \begin{equation}\label{eq:Psi_representation}
    \psi'(u) = \mathrm{i}\beta - u\sigma^2 +\mathrm{i} \int_{\mathbb{R}}\left( \frac{1}{x-\mathrm{i}u}-\frac{x}{1+x^2} \right)U (dx).
    \end{equation}
    Since $\psi(u)\to 0$ when $u\downarrow 0$, we must have that $U (\lbrace 0\rbrace)=0$. Indeed, for $r> 0$ it holds that
    \[
    -\Re (\psi'(r)) = \sigma^2 r + \int_{\mathbb{R}}\frac{r}{x^2 +r^2}U (dx). 
    \]
    Integrating this from 0 to $r$ and using $\psi(0)=0$ gives that
    \[
    -\log (\lvert \varphi(r)\rvert) =\frac{\sigma^2 r^2}{2}+\frac{1}{2}\int_{\mathbb{R}}\log\left( 1+\frac{r^2}{x^2} \right)U (dx).
    \]
    The left hand side is finite since $\varphi(r)\neq 0$ and thus $U(\lbrace 0 \rbrace)=0$. It also follows that, since $\log(1+\tfrac{r^2}{x^2})\asymp \lvert \log(\lvert x\rvert)\rvert$ as $x\to 0$, 
    \begin{equation}\label{eq:U_log_near_zero}
    \int_{\lbrace \lvert x\rvert\leq 1\rbrace}\lvert\log(\lvert x\rvert)\rvert U (dx)<\infty.
    \end{equation}
    Together with the integrability condition on $U$, we find that
    $
    \int_{\lbrace \lvert x\rvert\geq 1\rbrace}x^{-2}U (dx) < \infty.
    $
    Integrating \eqref{eq:Psi_representation} to obtain $\psi$ and absorbing the change of truncation into $\beta$ to  yield a constant term $\gamma\in\mathbb{R}$, we obtain the Thorin representation \eqref{eq:thorin_class_def} as (up to a truncation term that can be absorbed in $\gamma$)
    \begin{equation}\label{eq:log_phi_thorin_repr}
    \log (\varphi(u)) = \mathrm{i}\gamma u - \frac{\sigma^2 u^2}{2}+\int_{\mathbb{R}\setminus \lbrace 0 \rbrace} \left[\log\left( \frac{x}{x-\mathrm{i}u}\right)-\frac{\mathrm{iu}}{x}\mathbf{1}_{\lbrace \lvert x \rvert\geq 1\rbrace}(x)\right] U (dx).
    \end{equation}
    This requires an application of Fubini's theorem, which is justified by \eqref{eq:U_log_near_zero} and items 1 and 2 of Lemma~\ref{lem:U_integrability}. Absorbing the change of truncation into a constant is justified by item 3 in Lemma~\ref{lem:U_integrability}. This proves \eqref{eq:log_phi_thorin_repr} for $u>0$, but the extension to $u<0$ is simple. The right hand side of \eqref{eq:log_phi_thorin_repr} depends on $u$ only through $\mathrm{i}u$ and the principal logarithm commutes with conjugation away from the negative real axis. As such, exponentiating \eqref{eq:log_phi_thorin_repr}, evaluating at $-u$, and conjugating, yields \eqref{eq:log_phi_thorin_repr} for $u<0$ and hence all $u$.
\end{proof}

\begin{lemma}\label{lem:HFL_cf_is_Thorin}
    Let $\varphi$ be a characteristic function with the HBB property \eqref{eq:HFL_property} and such that its restriction to $(0,\infty)$ extends analytically to $\mathcal{H}=\lbrace u\in \mathbb{C}:\Re(u)> 0\rbrace$. Then $\varphi$ is the Fourier transform of a law in $\mathcal{T}(\mathbb{R})$.
\end{lemma}
\begin{proof}
    Let $\mu$ be the corresponding Bochner-Bernstein measure so that \eqref{eq:HFL_property} holds for $\varphi$. We first extend $\varphi$ to $-\mathcal{H}$ by
    \(
        \varphi(z):=\overline{\varphi(-\overline z)},
    \)
    for $z\in-\mathcal{H}$. This extension is holomorphic on $-\mathcal{H}$, agrees with the original characteristic function on \((-\infty,0)\), and satisfies
    $
        \varphi(-\overline z)=\overline{\varphi(z)},
    $ for $z\in\mathcal{H}$. 
    Fix $c>0$ and define the finite positive measure
    $
        \eta_c(dx)=\int_{[0,\infty)}e^{-c^{2}y}\mu(dx,dy).
    $
    Following Remark~\ref{rem:HFL_property}, let $p\in\mathbb{C}$ and set
    \[
        S_c^{+}(p)=\frac{p+\sqrt{p^{2}+4c^{2}}}{2},\qquad S_c^{-}(p)= \frac{-p+\sqrt{p^{2}+4c^{2}}}{2},
    \]
    where the principal square root is used. For real $p$, both $S_c^{+}(p)$ and $S_c^{-}(p)$ are positive, and
    \[
        S_c^{+}(p)-S_c^{-}(p)=p,\qquad S_c^{+}(p)S_c^{-}(p)=c^{2}.
    \]
    Consequently, \eqref{eq:HFL_property} gives
    \begin{equation}\label{eq:fixed-c-real}
        \widehat{\eta_c}(p)=\int_{\mathbb{R}}e^{\mathrm{i}px}\eta_c(dx)=\varphi\bigl(S_c^{+}(p)\bigr)\varphi\bigl(-S_c^{-}(p)\bigr),\quad p\in\mathbb{R}.
    \end{equation}
    The right-hand side of \eqref{eq:fixed-c-real} is holomorphic on the strip $\Sigma_c=\{p\in\mathbb{C}:|\operatorname{Im}p|<2c\}$. Indeed, \(p^{2}+4c^{2}\notin(-\infty,0]\) for \(p\in\Sigma_c\), so the
    principal square root is holomorphic there. Moreover,
    \[
        S_c^{+}(p),S_c^{-}(p)\in \mathcal{H},
        \quad p\in\Sigma_c.
    \]
    To see this, note that neither $S_c^{+}$ nor $S_c^{-}$ can cross the imaginary axis inside $\Sigma_c$. If, for example,
    $S_c^{+}(p)=\mathrm{i}r$, then
    $    
        S_c^{-}(p)=\tfrac{c^{2}}{S_c^{+}(p)}=-\tfrac{\mathrm{i}c^{2}}{r},
    $
    and hence
    \(
        \lvert\Im (p)\rvert=\lvert r+\frac{c^{2}}{r}\rvert\geq 2c,
    \)
    which is in contradiction to $p\in\Sigma_c$. The function $\widehat{\eta_c}$ admits a holomorphic extension to $\Sigma_c$ given by the right hand side of \eqref{eq:fixed-c-real}, which we denote by $A_c$ and note that for real $p\in\mathbb{R}$ $A_c(p)=\widehat{\eta}_c(p)$. In particular, the real Fourier transform $\widehat{\eta_c}$ admits a holomorphic extension to the disc $\lbrace \lvert p\rvert < 2c\rbrace$. Hence $\widehat{\eta}_c(p)$ is infinitely differentiable at 0 and we claim that $m_{2k}:=\int x^{2k}\eta_c(dx)<\infty$ and that $\widehat{\eta}_c^{(2k)}(0)=(-1)^k m_{2k}$ for every $k\geq 0$. This is clear for $k=0$ and if it holds for $k-1$, then $\widehat{\eta}_c^{(2k-2)}(p)=(-1)^{k-1}\int x^{2k-2}e^{\mathrm{i}px}\eta_c(dx)$ for all real $p$ due to dominated convergence. Hence we find for $h\neq 0$ that 
    \[
    (-1)^k \int x^{2k-2}\frac{2(1-\cos(hx))}{h^2}\eta_c(dx) = \frac{\widehat{\eta}_c^{(2k-2)}(h)-2\widehat{\eta}_c^{(2k-2)}(0)+\widehat{\eta}_c^{(2k-2)}(-h)}{h^2},
    \]
    where the right hand side tends to $\widehat{\eta}_c^{(2k)}(0)$ as $h\to 0$. As $2(1-\cos(hx))/h^2 \to x^2$ as $h\to 0$ and the integrands are non-negative, we may apply Fatou's lemma to find that $m_{2k}\leq (-1)^k \widehat{\eta}_c^{(2k)}(0)<\infty$. Dominated convergence then upgrades the inequality to an equality. We then find that, for $0\leq r<2c$,
    \[
    \int_\mathbb{R}\cosh(rx)\eta_c(dx)=\sum_{k\geq 0}\frac{r^{2k}}{(2k)!}m_{2k}=\sum_{k\geq 0}\frac{r^{2k}}{(2k)!}\lvert A_c^{(2k)}(0)\rvert <\infty,
    \]
    where we have used that the Taylor series of $A_c$ is absolutely convergent at $\lvert p\rvert =r$. Since, $e^{r\lvert x\rvert}\leq 2\cosh(rx)$, we find that
    \begin{equation}\label{eq:exponential_moment_eta}
        \int_{\mathbb{R}}e^{\alpha |x|}\eta_c(dx)<\infty,\quad 0\leq \alpha<2c,
    \end{equation}
    for every $c>0$. We can therefore define, for every $p\in\Sigma_c$, the function $F_c(p)=\int_\mathbb{R}e^{\mathrm{i}px}\eta_c(dx)$, where the integral is absolutely convergent since $\lvert e^{\mathrm{i}px}\rvert\leq e^{\lvert \Im (p)\rvert \lvert x\rvert}$ and $\lvert \Im(p)\rvert<2c$ on $\Sigma_c$ by definition. We claim that $F_c=A_c$ on $\Sigma_c$, which means that
    \begin{equation}\label{eq:proof1_hfl}
    \int_{\mathbb{R}\times [0,\infty)}e^{\mathrm{i}px-c^2y}\mu(dx,dy) = \varphi(S^+_c(p))\varphi(-S_c^-(p)).
    \end{equation}
    If this claim holds, then we may set $z=a+\mathrm{i}b\in\mathcal{H}$ (so $a>0$) and set $c=\sqrt{a^2+b^2}$ and $p=2\mathrm{i}b$. It then holds that $\lvert \Im (p)\rvert =2\lvert b\rvert <2c$ and that
    $
    \sqrt{p^2 +4c^2}=\sqrt{-4b^2+4(a^2+b^2)}=2a.
    $
    Therefore, $S_c^+(p)=a+\mathrm{ib}=z$ and $S_c^-(p)=a-\mathrm{i}b=\overline{z}$. Evaluating \eqref{eq:proof1_hfl} at this $p$ yields
    \begin{equation}\label{eq:phi_modulus_representation}
    \lvert \varphi(a+\mathrm{i}b)\rvert^2 = \int_{\mathbb{R}\times [0,\infty)}e^{-2bx-(a^2+b^2)y}\mu (dx,dy) > 0,
    \end{equation}
    such that $\varphi$ has no zeros in $\mathcal{H}$. Then, fixing $b\in\mathbb{R}$ and differentiating \eqref{eq:phi_modulus_representation} with respect to $a>0$ yields
    \begin{equation}\label{eq:diff_modulus}
    2\Re \left( \varphi'(a+\mathrm{i}b)\overline{\varphi(a+\mathrm{i}b} )\right) = -2a\int_{\mathbb{R}\times [0,\infty)}ye^{-2bx-(a^2+b^2)y}\mu(dx,dy),
    \end{equation}
    where differentiation under the integral sign is justified below. Dividing by $2\lvert \varphi (a+\mathrm{i}b)\rvert^2 >0$ then yields 
    \[
        \Re\left(\frac{\varphi'(a+\mathrm{i}b)}{\varphi(a+\mathrm{i}b)}\right)=-a
        \frac{\int ye^{-2bx-(a^{2}+b^{2})y}\mu(dx,dy)
        }{\int e^{-2bx-(a^{2}+b^{2})y}\mu(dx,dy)} \leq 0.
    \]
    Hence $\varphi$ is non-vanishing on $\mathcal{H}$ and satisfies $\Re (\tfrac{\varphi'(z)}{\varphi(z)})\leq 0$ for $z\in \mathcal{H}$ and the conclusion follows from Lemma~\ref{lem:pick_theorem}.
    
    It remains to argue that $F_c=A_c$ on $\Sigma_c$ and to justify the differentiation under the integral sign to obtain \eqref{eq:diff_modulus}. First, we show that $F_c$ is holomorphic on $\Sigma_c$. Fix $p_0\in \Sigma_c$ and constants
    \[
    \lvert \Im (p_0)\rvert<\gamma<\alpha <2c.
    \]
    There is a neighbourhood $U$ of $p_0$ such that $\lvert\Im (p)\rvert \leq \gamma$ for all $p\in U$ and for such $p$ it holds that
    \[
        \lvert x e^{\mathrm{i}px}\rvert \leq \lvert x\rvert e^{\gamma \lvert x\rvert} =\left(\lvert x \rvert e^{-(\alpha-\gamma)\lvert x\rvert}\right)e^{\alpha\lvert x \rvert} \leq C_{\alpha,\gamma}e^{\alpha\lvert x \rvert},
    \]
    for some constant $C_{\alpha,\gamma}>0$. The final expression is integrable with respect to $\eta_c$ due to \eqref{eq:exponential_moment_eta}. Differentiation under the integral sign is therefore justified locally around $p_0$ and
    \[
        F_c'(p)=\mathrm{i}\int_{\mathbb R}x e^{\mathrm{i}px}\eta_c(dx).
    \]
    Therefore, $F_c$ is holomorphic on $\Sigma_c$ and so $F_c$ and $A_c$ are holomorphic on the connected domain $\Sigma_c$. For real $p$, equation \eqref{eq:fixed-c-real} gives that $F_c(p)=A_c(p)$. Since the real line has accumulation points in $\Sigma_c$, the identity theorem for analytic functions yields that $F_c(p)=A_c(p)$ for all $p\in \Sigma_c$, as claimed. 
    
    To justify differentiation under the integral sign in \eqref{eq:diff_modulus}, we fix $a_0>0$ while keeping $b\in\mathbb{R}$ fixed and let $d$ be such that $\lvert b\rvert < d<\sqrt{a_0^2+b^2}$. By continuity, there is a neighbourhood $I$ of $a_0$ and a constant $\delta>0$ such that $a^2+b^2 \geq d^2 +\delta$ for all $a\in I$. We may also choose $M<\infty$ such that $a\leq M$for every $a\in I$, possibly by first shrinking $I$. Consequently,
    \[
    \begin{aligned}
        \left\lvert \frac{\partial}{\partial a} e^{-2bx-(a^2+b^2)y}\right\rvert
        &=2aye^{-2bx-(a^2+b^2)y}\\
        &\leq2Mye^{2|b||x|-(d^2+\delta)y}\\
        &=2M\bigl(y e^{-\delta y}\bigr)e^{2|b||x|-d^2y}\\
        &\leq C_{\delta,M}e^{2|b||x|-d^2y}.
    \end{aligned}
    \]
    The dominating function on the right hand side is integrable since 
    \[
        \int_{\mathbb{R}\times[0,\infty)}e^{2\lvert b\rvert\lvert x\rvert-d^2y}\mu(dx,dy)= \int_{\mathbb R}e^{2\lvert b \rvert \lvert x \rvert}\eta_d(dx)<\infty,
    \]
    where the final inequality follows from \eqref{eq:exponential_moment_eta} with $c=d$. Dominated convergence then justifies \eqref{eq:diff_modulus}. 
\end{proof}

\begin{lemma}\label{lem:gamma_multiplier_analytic}
    Let $X$ be an $\mathbb{R}$-valued random variable and let $\Gamma_n\sim \operatorname{Gamma}(n,n)$ be independent of $X$. Then the characteristic function of $\Gamma_n X$ extends analytically to $\mathcal{H}=\lbrace u\in \mathbb{C}:\Re (u)>0\rbrace$.
\end{lemma}
\begin{proof}
    The characteristic function of $\Gamma_nX$ is
    \[
        \varphi_{\Gamma_nX}(u) = \int_0^\infty \varphi_{X}(ug)\frac{n^n}{\Gamma (n)}g^{n-1}e^{-ng}dg = \frac{n^n}{\Gamma (n)}u^{-n}\int_0^\infty \varphi_X (w)w^{n-1}e^{-\tfrac{n}{u}w}dw.
    \]
    Note that $\Re (n/u)=n\Re (u/\lvert u\rvert^2)>0$ and $\lvert \varphi_X\rvert\leq 1$, so the integral on the right hand side converges absolutely and locally uniformly. Hence $\varphi_{\Gamma_nX}$ extends analytically to all of $\mathcal{H}$.
\end{proof}

\subsection{Proof of Theorem~\ref{thm:HCM_multiplication_closure}}
Let $Y'$ be an independent copy of $Y$ and let $(X,X')$ be the coupling of Proposition~\ref{prop:Thorin_coupling}. Equation~5.1.1  of \textcite{Bondesson92} shows that the HCM density $f$ of $Y$ satisfies for every $\rho >0$
    \[
    f(\rho v)f(\rho/v) = \int_{[0,\infty)}e^{-\lambda\rho (v+v^{-1})}\zeta(d\lambda;\rho),\quad v>0,
    \]
    where $\zeta$ is a positive measure. Let $\varphi_Z$ denote the characteristic function of a random variable $Z$ and note that the independence of $(X,X')$ and $(Y,Y')$ allows for writing
    \begin{equation}\label{eq:cf_XY_exp}
    \begin{aligned}
    \varphi_{YX}(u)\varphi_{YX}(-w)&=\mathbb{E}\left[ \varphi_X(uY)\varphi_X(-wY') \right] \\
    &= \mathbb{E}\left[ e^{\mathrm{i}(uY-wY')X-uwYY'X'} \right].
    \end{aligned}
    \end{equation}
    Since $(X,X'),(Y,Y')$ are independent, we can realize $(X,X')$ and $(Y,Y')$ on separate probability spaces and obtain
    \begin{equation}\label{eq:FL_exp1}
    \mathbb{E}\left[ e^{\mathrm{i}(uY-wY')X-uwYY'X'} \right] = \mathbb{E}\left[ \int_{0}^{\infty}\int_0^\infty e^{\mathrm{i}(uy-wy')X-uwyy'X'} f(y)f(y')\,dy\,dy' \right].
    \end{equation}
    The hyperbolic change of variables $\rho = \sqrt{yy'}$ and $v=\sqrt{y/y'}$ yields that $y'=\rho/v$ and $y=\rho v$ and substituting into the right hand side of \eqref{eq:FL_exp1} yields
    \[
    \begin{aligned}
    \varphi_{YX}(u)\varphi_{YX}(-w) =\mathbb{E}\left[ \int_0^\infty\int_0^\infty 2\frac{\rho}{v} h(v)e^{-\rho^2uwX'}f(\rho v)f(\rho/v)\, d\rho \, dv \right],
    \end{aligned}
    \]
    where $h(v)=e^{\mathrm{i}\sqrt{uw}\rho X\left[v\sqrt{u/w}-1/(v\sqrt{u/w})\right]}$. Using the HCM property of $f$ we find
    \[
    \varphi_{YX}(u)\varphi_{YX}(-w) =\mathbb{E}\left[ \int_0^\infty\int_0^\infty 2\frac{\rho}{v} h(v)e^{-\rho^2uwX'}\int_{[0,\infty)}e^{-\lambda\rho (v+v^{-1})}\zeta(d\lambda;\rho)\, d\rho \, dv \right].
    \]
    Note that $\rho\mapsto \zeta (d\lambda;\rho)$ is measurable since $\widetilde{\zeta}(d\lambda;\rho)=e^{-2\lambda\rho}\zeta(d\lambda,\rho)$ is a finite measure with mass $f(\rho)^2$ and $\int e^{-\lambda s}\widetilde{\zeta}(d\lambda;\rho)=f(\rho v)f(\rho/v)$ with $v+1/v=2+s/\rho$ is Borel in $\rho$ for each $s\geq 0$; exponentials generate the Borel sets of $[0,\infty)$ so $\widetilde{\zeta}(A;\cdot)$ is Borel for every $A$ and so is $\zeta(A;\cdot)$.
    Furthermore, it must hold that $\zeta(\lbrace 0\rbrace ; \rho)=0$ for almost every $\rho$ since an atom at $0$ would contribute with $\zeta(\lbrace 0\rbrace;\rho)\int \tfrac{1}{v}dv=\infty$. We want to first consider the integral with respect to $v$, which requires interchanging the order of integration. To this end, note that the integrand is bounded in modulus by $\tfrac{2\rho}{v}e^{-\lambda \rho (v+v^{-1})}$ whose full integral is simply $\iint f(y)f(y')dydy' = 1$ (seen by undoing the hyperbolic substitution). Fubini's theorem thus applies and we may swap the order of integration and treat the integral with respect to $v$. To this end, we introduce the quantities
    \[
    A=\lambda\rho - \mathrm{i}u\rho X, \quad B=\lambda\rho +\mathrm{i}w\rho X.
    \]
    Define $\Delta=\lbrace \Re (A)>0,\Re(B)>0,\Re(AB)>0\rbrace$ and note that the following identity holds $\int_0^\infty x^{\nu -1}e^{-\beta /x-\gamma x}dx = 2(\beta/\gamma)^{\nu /2}K_\nu(2\sqrt{\beta\gamma})$ whenever $\Re(\beta),\Re(\gamma)>0$ and $K_\nu$ denotes the modified Bessel function of the second kind. We therefore find that, on $\Delta$, it holds that
    \begin{equation}\label{eq:AB_integral}
    \int_0^\infty e^{-Av-B/v}\frac{1}{v}dv=\int_0^\infty e^{-\rho'-AB/\rho'}\frac{1}{\rho'}d\rho'.
    \end{equation}
    This representation is convenient, since
    \[
    AB= \lambda^2\rho^2 + q\rho^2X^2 - \mathrm{i}p\lambda \rho^2 X,
    \]
    where $q=uw$ and $p=u-w$. Inserting this into \eqref{eq:AB_integral} we obtain
    \begin{equation}\label{eq:v_integral}
    \int_0^\infty h(v)e^{-\lambda \rho (v+v^{-1})}\frac{1}{v}dv = \int_0^\infty e^{-\rho'-\lambda^2\rho^2/\rho'}e^{-q\frac{\rho^2 X^2}{\rho'}}e^{\mathrm{i}p\frac{\lambda\rho^2 X}{\rho'}}\frac{1}{\rho'}\, d\rho'.
    \end{equation}
    I.e., the $v$ integral is a positive weight times a Laplace kernel in $q$ at the non-negative point $\frac{\rho^2 X^2}{\rho'}$ times a Fourier kernel in $p$ at the point $\frac{\lambda\rho^2 X}{\rho'}$. Note that this representation only holds on $\Delta$, which requires $\Re (AB)>0$. In our case this condition becomes $\lambda^2\rho^2+q\rho^2 X^2>0$, which in fact explains why the conjugate $\varphi(-w)$ and not $\varphi(w)$ is the right object.

    Now define the measure $M(d\rho,d\lambda,d\rho',dx,dx')$ as
    \[
    M(d\rho,d\lambda,d\rho',dx,dx') = \frac{2\rho}{\rho'}e^{-\rho'-\lambda^2\rho^2/\rho'}\zeta(d\lambda;\rho)d\rho d\rho' \mathbb{P}\left((X,X')\in (dx,dx')\right).
    \]
    Note that $\mathbb{P}\left( (X,X')\in \mathbb{R}\times\mathbb{R}_+ \right)=1$ and that 
    \[
    \int_{0}^{\infty}e^{-\rho'-\lambda^2\rho^2/\rho'}\frac{1}{\rho'}\, d\rho' = \int_{0}^{\infty}e^{-\lambda\rho (v+v^{-1})}\frac{1}{v}dv.
    \]
    The remaining integral over $(\lambda,\rho)$ is thus also 1, since undoing the hyperbolic substitution used above yields
    \[
    \int_0^{\infty}\left(\int_0^{\infty}\left(\int_0^\infty e^{-\lambda\rho (v+1/v)}\frac{1}{v}dv\right) \zeta (d\lambda;\rho)\right)2\rho d\rho = \int_0^\infty\int_0^\infty f(y)f(y')dydy' = 1.
    \]
    Consequently, $M$ has total mass $1$ and we may define a probability measure $\mu$ on $\mathbb{R}\times [0,\infty)$ as the image of $M$ under the map
    \begin{equation}\label{eq:M_pushforward}
    (\rho,\lambda,\rho',x,x')\mapsto (z,z') := \left(\tfrac{\lambda \rho^2 x}{\rho'}, \, \rho^2 x'+\tfrac{\rho^2 x^2}{\rho'}\right).
    \end{equation}
    Note that it follows from \eqref{eq:v_integral} and the definition of $M$ that
    \[
    \varphi_{YX}(u)\varphi_{YX}(-w)=\int e^{\mathrm{i}(u-w)\frac{\lambda\rho^2 x}{\rho'}-uw\left(\rho^2 x' +\frac{\rho^2 x^2}{\rho'}\right)}M (d\rho,d\lambda, d\rho',dx,dx'), \quad u,w>0.
    \]
    Pushing forward the measure $M$ via \eqref{eq:M_pushforward}, the change-of-variables formula yields that
    \[
    \varphi_{YX}(u)\varphi_{YX}(-w)=\int_\mathbb{R}\int_0^\infty e^{\mathrm{i}(u-w)z-uwz'}\mu (dz,dz'), \quad u,w>0,
    \]
    which is exactly the Bochner-Bernstein property in the coordinates $(z,z')$ and thus we conclude that $XY\in\mathcal{H}_{BB}$. To see that $X/Y\in \mathcal{H}_{BB}$, note that Proposition~4 of \textcite{Bondesson2015} proves that $1/Y\in\mathcal{H}_{CM}$, so $X/Y\in\mathcal{H}_{BB}$ follows from the closure with respect to products.

\subsection{Proof of Theorem~\ref{thm:Thorin_HFL_equivalent}} 
    Set $X\sim \operatorname{Gamma}(\alpha,\beta)$ and note that for $c\in\mathbb{C}$ with $\Re (c)>-\beta$ it holds that $\mathbb{E}[e^{-cX}]=(1+c/\beta)^{-\alpha}$ where the complex root is taken as the principal branch. Set $c=\tfrac{uw}{\beta}-\mathrm{i}(u-w)$ for $u,w>0$ such that $\Re (c)>-\beta$ and note that $\arg(1-\tfrac{\mathrm{iu}}{\beta})\in (-\pi/2,0)$ and $\arg (1+\tfrac{\mathrm{i}w}{\beta})\in (0,\pi/2)$, such that the principal powers multiply and it therefore holds that
    \[
    \mathbb{E}\left[ e^{\mathrm{i}(u-w)X-\tfrac{uw}{\beta}X} \right] = \left( 1+\frac{uw}{\beta^2}-\frac{\mathrm{i}(u-w)}{\beta} \right)^{-\alpha}=\left(1-\frac{\mathrm{iu}}{\beta}\right)^{-\alpha}\left( 1+\frac{\mathrm{i}w}{\beta}\right)^{-\alpha}=\varphi_X(u)\varphi_X(-w).
    \]
    Thus the coupling of Proposition~\ref{prop:Thorin_coupling} is $(X,X/\beta)$. Swapping the variables $u$ and $w$ yields the same result for the random variable $\widetilde{X}=-X$ and thus it holds that any two-sided gamma law has the HBB property. Letting $Y\in \mathcal{T}(\mathbb{R})$, $Y$ is a weak limit of two-sided gamma convolutions. Convolution closure and closure under weak limits of HBB laws cf. Lemma~\ref{lem:HFL_closure} yields that $Y\in \mathcal{H}_{BB}$ and hence every law in Thorin's class has the HBB property.

    Conversely, let $X\in \mathcal{H}_{BB}$ and set $Z_n=\Gamma_n X$ with $\Gamma_n\sim \operatorname{Gamma}(n,n)$ independent of $X$. Since $\Gamma_n\in \mathcal{H}_{CM}$, it follows from Theorem~\ref{thm:HCM_multiplication_closure} that $Z_n\in \mathcal{H}_{BB}$. The characteristic function $\varphi_{Z_n}$ is analytic on $\mathcal{H}$ by Lemma~\ref{lem:gamma_multiplier_analytic} and therefore $Z_n\in\mathcal{T}(\mathbb{R})$ by Lemma~\ref{lem:HFL_cf_is_Thorin}. Since $\Gamma (n,n)\xrightarrow[n\to\infty]{\mathbb{P}}1$, it follows that $Z_n\xrightarrow[n\to\infty]{d}X$ and since $\mathcal{T}(\mathbb{R})$ is closed under weak limits, we find that $X\in\mathcal{T}(\mathbb{R})$.   

\section*{Funding}
This work was supported by the Center of Research in Energy:
Economics and Markets and The Danish Council of Independent Research under DFF grant 10.46540/5247-00005B

\section*{Disclosure statement}
The author declares that there are no conflicts of interest.

\printbibliography

\end{document}